\documentclass{amsart}
 \usepackage{amsmath,amsthm,amssymb,amscd, amstext, amsopn, amsxtra}

\newtheorem{thm}{Theorem}[section]
\newtheorem*{thmA}{Theorem}
\newtheorem{lemma}{Lemma}
\newtheorem{cor}[thm]{Corollary}

\newtheorem{claim}{Claim}
\newtheorem{prop}{Proposition}
\theoremstyle{remark}
\newtheorem{defn}{Definition}
\newtheorem{rul}{Rule}
\newtheorem{rem}{Remark}

\newcommand{\Z}{{\mathbb Z}}

\newcommand{\1}{{\bf 1}} 

\newcommand{\B}{B}

\newcommand{\ttau}{s}
\newcommand{\tensor}{\otimes}
\newcommand{\ot}{\otimes^*}

\DeclareMathOperator{\Rep}{Rep}

\DeclareMathOperator{\Res}{Res}
\newcommand{\Hom}{{\mathop{\mathrm{Hom}}\nolimits}}

\newcommand{\Ind}{{\mathop{\mathrm{Ind}}\nolimits}}
\newcommand{\aff}{{\mathop{\mathrm{aff}}\nolimits}}

\newcommand{\Indhat}{{\mathop{\mathrm{\widehat Ind}}\nolimits}}

\newcommand{\fin}{{\mathop{\mathrm{fin}}\nolimits}}
\newcommand{\isom}{\simeq}
\newcommand{\surj}{\twoheadrightarrow}
\newcommand{\inj}{\hookrightarrow}

\newcommand{\rev}{{\mathop{\mathrm{rev}}\nolimits}}

\newcommand{\ee}{e}
\newcommand{\h}{h}

\newcommand{\e}[1]{\ee_{#1}}
\newcommand{\ei}{\e{i}} 
\newcommand{\ff}{f}
\newcommand{\et}[1]{{\widetilde{\ee}_{#1}}} 
\newcommand{\eti}{\et{i}} 
\newcommand{\ft}[1]{\widetilde{\ff}_{#1}}
\newcommand{\fti}{\ft{i}}
\newcommand{\epsi}{\varepsilon_i}
\newcommand{\epsj}[1]{\varepsilon_{#1}}
 
\DeclareMathOperator{\soc}{soc}
\DeclareMathOperator{\cosoc}{cosoc}
\DeclareMathOperator{\ch}{ch}

\newcommand{\RqHaffn}{\Rep_q \Haffn}
\newcommand{\Rq}{\Rep_q}

\newcommand{\RqHaff}{\Rep_q H_n^{\aff}}

\newcommand{\RLami}{\Rep H_n^{\Lambda_i}}

\newcommand{\RHlam}{\Rep H_n^{\lambda}}
\newcommand{\RHlamn}{\Rep H_n^{\lambda}}

\newcommand{\Haffn}{\Haff{n}}
\newcommand{\Haff}[1]{H_{#1}}
\newcommand{\Hlam}[1]{H_{#1}^{\lambda}}
\newcommand{\Hlamn}{\Hlam{n}}

\newcommand{\gl}{{\mathfrak gl}_\infty}
\newcommand{\Rx}{R[X_1^{\pm 1}, \ldots, X_n^{\pm 1}] }

\newcommand{\q}[2]{{ q^{#1}  \cdots q^{#2} }}
\newcommand{\qq}[2]{{ q^{#1} q^{{#1} +1} \cdots q^{#2} }}
\newcommand{\del}[2]{\triangle_{(#1,#2)}}
\newcommand{\Del}{{\mathbf \Delta}}

\newcommand{\bx}{\boxtimes} 

\newcommand{\ehat}[1]{\widehat{\ee}_{#1}}
\newcommand{\epshat}[1]{\widehat{\varepsilon}_{#1}}
\newcommand{\phati}{\widehat{\varphi}_{i}}
\newcommand{\epshati}{\epshat{i}}
\newcommand{\ehati}{\ehat{i}}

\newcommand{\ethat}[1]{{\widehat{\tilde{\ee}}}_{#1}}
\newcommand{\ethati}{\ethat{i}}

\newcommand{\qi}{q^i}

\newcommand{\aaa}{a}

\newcommand{\WP}[1]{W^{#1}_{\min}}

\newcommand{\oK}{\bigoplus_n K(\RqHaffn)}
\newcommand{\resn}{\Res_{H_{n-1} \otimes H_1}^{H_n}} 

\newcommand{\ve}{\varepsilon}
\def\std{left }
\def\lex{right }

\newcommand{\xDel}{\Del}
\newcommand{\sDel}{\underline{\Del}}
\newcommand{\sGamma}{\underline{\Gamma}}
\newcommand{\Mg}[1]{M_{#1}}
\newcommand{\Md}{\Mg{\Del}}
\newcommand{\Delsup}[1]{\Del^{(#1)}}
\newcommand{\Delsupj}{\Delsup{j}}
\newcommand{\Delbar}{\Del^{-}}
\newcommand{\Qdel}{Q}
\newcommand{\E}[1]{{\mathbf {\et{#1}} }}
\newcommand{\Ei}{\E{i}}
\newcommand{\F}[1]{{\mathbf {\ft{#1}} }}
\newcommand{\Fi}{\F{i}}

\newcommand{\Ehat}[1]{{\mathbf {\ethat{#1}} }}
\newcommand{\Ehati}{\Ehat{i}}
\newcommand{\0}{{\mathbf 0}}
\newcommand{\Nm}[2]{{N_{{#1},{#2}}}}
\newcommand{\Nmb}[1]{N_{{#1}}}
\newcommand{\Nmui}{\Nm{\mu}{i}}
\newcommand{\Nmulam}{\Nmb{\mubarlam}}
\newcommand{\mubar}[2]{\underline{#1}, {#2}}
\newcommand{\mubari}{\mubar{\mu}{\lambda}}
\newcommand{\mubarlam}{\mubar{\mu}{\lambda}}
\newcommand{\mubarik}{\underline{\mu}}
\newcommand{\nubarik}{\underline{\nu}}
\newcommand{\mubaralpha}{\mubar{\mu}{\alpha}}
\newcommand{\nubarbeta}{\mubar{\nu}{\beta}}
\newcommand{\concat}[2]{{#1}\cdot{#2}}
\newcommand{\vpi}{\varphi_i}

\newcommand{\jboxj}{\framebox{\makebox[\totalheight]{j}}}

\def\shuffle{\,\raise 1pt\hbox{$\scriptscriptstyle\cup{\mskip
               -4mu}\cup$}\,}

\newlength\cellsize 
\savebox2{%
\begin{picture}(18,18)
\put(0,0){\line(1,0){18}}
\put(0,0){\line(0,1){18}}
\put(18,0){\line(0,1){18}}
\put(0,18){\line(1,0){18}}
\end{picture}}

\newcommand\cellify[1]{\def\thearg{#1}\def\nothing{}%
\ifx\thearg\nothing
\vrule width0pt height\cellsize depth0pt\else
\hbox to 0pt{\usebox2\hss}\fi%
\vbox to 18\unitlength{
\vss
\hbox to 18\unitlength{\hss$#1$\hss}
\vss}}

\newcommand\tableau[1]{\vtop{\let\\=\cr
\setlength\baselineskip{-16000pt}
\setlength\lineskiplimit{16000pt}
\setlength\lineskip{0pt}
\halign{&\cellify{##}\cr#1\crcr}}}

\savebox3{%
\begin{picture}(15,15)
\put(0,0){\line(1,0){15}}
\put(0,0){\line(0,1){15}}
\put(15,0){\line(0,1){15}}
\put(0,15){\line(1,0){15}}
\end{picture}}
\newcommand\expath[1]{%
\hbox to 0pt{\usebox3\hss}%
\vbox to 15\unitlength{
\vss
\hbox to 15\unitlength{\hss$#1$\hss}
\vss}}
 
\newcommand{\omitt}[1]{}

 \title[Parameterizing Hecke algebra modules]{Parameterizing
Hecke algebra modules:
Bernstein-Zelevinsky multisegments,
Kleshchev multipartitions, and crystal graphs
}

 \author{M. Vazirani}
 \date{\today}

\begin{document}
 \maketitle
 \tableofcontents
 \section{Introduction} \label{sec-intro}
This paper provides a combinatorial dictionary between
three sets of objects:  multisegments, multipartitions, and
the irreducible modules of the affine Hecke algebra
$\Haffn$ (for generic $q$).
The dictionary is dictated by Grojnowski's 
Theorem 14.3,  \cite{G}, (repeated here as
Theorems \ref{thm-graph} and \ref{thm-graphinfty})
 in which he constructs the crystal
graph with nodes given by the irreducible modules of $\Haffn$
and proves that graph is $\B(\infty)$ (Kashiwara's crystal
graph associated to $U(\eta^-)$).
He also shows the subgraph with nodes the irreducible modules
of $\Hlamn$ is $\B(\lambda)$, the crystal graph associated to
the irreducible highest weight representation of $\gl$ with highest weight
$\lambda$.

One can also construct purely combinatorial crystal graphs
whose nodes are Bernstein-Zelevinsky multisegments,
or are Kleshchev multipartitions, and these crystal graphs
are {\it abstractly\/} isomorphic to $\B(\infty)$ and $\B(\lambda)$, 
respectively.
In this paper we give explicit isomorphisms.
In particular,
we compute the action of the crystal operator $\eti$ on an
irreducible module both in terms of its parameterization by 
multisegments (theorem \ref{thm-Ei}, rule \ref{defn-Ei})
and by multipartitions  (theorem \ref{thm-part2}, rule \ref{defn-p}).

The results below contain new representation-theoretic content.
This is because, following \cite{BZ}, we directly define a map
from multisegments to irreducible modules
$$\Del \mapsto \cosoc \Ind \Del$$
(the notation is explained in section \ref{sec-def1} below).
We 
then prove that this map intertwines the action of $\Ei$ on
multisegments (rule \ref{defn-Ei} of section \ref{sec-defnmulti})
with the action of $\eti$ on modules (section \ref{sec-ei}).
Similarly, we define a map from 
$\lambda$-colored Kleshchev multipartitions
(section \ref{sec-defnpart})
to irreducible modules
$$
(\mu^{(1)}, \ldots, \mu^{(r)}) \mapsto 
\cosoc \Ind \Nm{\mu^{(1)}}{i_1} \bx \cdots \bx \Nm{\mu^{(r)}}{i_r}
$$
and prove that this map  also intertwines the action of $\Ei$.
This is new.

Grojnowski's  Theorem 14.3  suffices to parameterize the modules
in terms of paths to the highest weight node, i.e.~ by
a sequence of operators $\cosoc 
\Ind$. The advantage of the approach of this paper is that we
only have to take $\cosoc$ once. 
However, it is not true in general that $\cosoc$ commutes with $\Ind$,
and this is the representation-theoretic difficulty we overcome.

Another
 byproduct of Theorem \ref{thm-Ei} below
is the determination of which multisegments
parameterize modules of the {\it cyclotomic} Hecke algebra
$\Hlamn$   (Theorem \ref{thm-Ehat}).
The theorems also partially explain why the rule for computing
$\eti$ mirrors the rule we know for that on a tensor product of crystal graphs.

The proofs given here do not rely on \cite{G} Theorem 14.3, but do use several
results from \cite{G}, \cite{GV} (and of course from \cite{BZ} and \cite{Z}).
One reason for this is that
the proofs are relatively straightforward.
Another is that  several steps in the proofs given here are
necessary ingredients for proving Theorems \ref{thm-part1} and \ref{thm-part2}.

There are related results in the literature in 
\cite{A2, AM,  G2}.
Both our proofs and construction are different from theirs,
and we comment on this at appropriate points in the text.

This introduction and section \ref{sec-crystals}
interpret the main theorems of this paper into the language
of crystal graphs.
However, the rest of the sections avoid that language
and simply deal with algebras and modules.
No background in crystal graphs is necessary to understand
the proofs, (but is helpful in appreciating them).

We begin with some necessary definitions and notation, then immediately
state the main results.
However, the reader may want to read the more extensive definitions
is section \ref{sec-def2} before section \ref{sec-main}.
  In sections \ref{sec-proofseg} and
\ref{sec-proofpartition} the proofs are given. In section \ref{sec-bzseg}
some results and proofs
that appear in \cite{BZ,Z} are given, recast in this notation,
and included for convenience to the reader.

We point out that these proofs require $q$ to be generic, whereas 
the theorems of \cite{G} hold for all $q \in R^\times$.

 \section{Background notations and results} \label{sec-def1}
In the following subsections, we define
three sets of objects which will be the nodes of three 
different crystal graphs.
These are the irreducible modules of $\Haffn$ in the
subcategory  $\Rq$,
multisegments,
and colored  multipartitions. We define the action of an operator $\eti$
on each of these three sets, i.e.~ we describe the edges of the crystal graphs.
In section \ref{sec-main} are the theorems showing the three crystal graphs
are isomorphic.

\subsection{$H_n^\fin$, $\Haffn$ and $H_n^\lambda$}
\label{sec-H}
Throughout the paper we fix an algebraically closed
 field $R$, and an invertible
element $q \in R$, 
such that $q^\ell = 1$ implies $\ell = 0$. 
In this situation $\{ q^i \}_{i \in \Z}$ is infinite 
and we say  $q$ is {\it generic\/}.

The {\it finite Hecke algebra\/}, $H_n^\fin$ is the $R$-algebra
with generators
$$T_1, \ldots , T_{n-1}$$
and relations
\begin{eqnarray}
\label{braid}
 \text{{\it braid relations\/}} & 
 T_i T_{i+1} T_i = T_{i+1} T_i T_{i+1}, \quad
T_i T_j = T_j T_i,\; | i-j | > 1 &
\\
\label{quad}
\text{{\it  quadratic relations\/}  } &
(T_i + 1)(T_i - q) = 0. &
\end{eqnarray}

The braid relations imply that if
$w = \ttau_{i_1} \cdots \ttau_{i_k}$ and $\ell (w) = k$,
then $T_{i_1} \cdots T_{i_k}$ depends only on $w \in S_n$.
It is denoted $T_w$, and the $\{T_w \mid  w \in S_n\}$ form a basis
of $H_n^\fin$ over $R$.

The  {\it affine Hecke algebra\/}
\cite{BZ}
%
$\Haffn$ (or $H_n^{\aff}$)
 is the $R$-algebra,
which as an $R$-module is isomorphic to
$$H_n^\fin \otimes_R \Rx.$$
The algebra structure is given by requiring that $H_n^\fin$ and
$\Rx$ are subalgebras, and that
\begin{equation}
T_i X_i T_i = q X_{i+1} .
\label{txt}
\end{equation}

Denote by $\RqHaffn$ the finite dimensional modules $M$ for $\Haffn$ 
such that the only eigenvalues of the $X_j$ on 
$M$ are powers of $q$.
Let $\Rq = \bigoplus_{n \ge 0} \RqHaffn$.

\begin{rem} \label{rem-qlines}
It follows from the computation in \cite{G}, section 6.2,
or is explained in \cite{V2}, or many other places, that 
to understand $\Rep \Haffn$ it is enough to understand $\RqHaffn$.
\end{rem}

 The {\it cyclotomic Hecke algebra\/} or {\it Ariki-Koike algebra\/}
\cite{AK}
  $H_n^\lambda$ is the
quotient $\Haffn / I_\lambda$,
where $I_\lambda$ is the ideal generated by
the polynomial in $X_1$:
$\prod (X_1 - q^i)^{m_i}$,
and $\lambda =  \sum m_i \Lambda_i$ is a weight of $\gl$
(where $\{ \Lambda_i \mid i \in \Z \}$ denote the fundamental weights).
We may also write $\lambda = \sum_{k=1}^r \Lambda_{i_k}$
and so $I_\lambda$ is generated by 
$\prod_{k=1}^r (X_1 - q^{i_k})$.

 \begin{rem}
\label{rem-cycl}
Observe that
any irreducible module in $\Rq$
is an irreducible
$\Hlamn$-module if we take $\lambda$ large enough.
Just let $\prod_{i=1}^r (X_1 - \qi)^{m_i}$ be the
characteristic polynomial of $X_1$ acting on the module.
Conversely, we identify any irreducible $\Hlamn$-module as 
an irreducible $\Haffn$-module (on which $(X_1 - \qi)^{m_i}$ vanishes).
\end{rem}

\subsection{$\ei$, $\eti$, and $\epsi$}
\label{sec-ei}

Given an irreducible $\Haffn$-module $M$, consider 
$\Res_{n-1,1}^n M = \resn M$.
As $X_n -  \qi \subseteq Z(H_{n-1}\otimes H_1)$,
left multiplication by $(X_n -  \qi)^m$ induces an endomorphism
of $\resn M$, and its kernel for
$m \gg 0$  (we can take $m \le \dim M$) is the  generalized eigenspace
of $X_n - \qi$.
Define $\ei M$
to be the restriction of that generalized eigenspace to $\Haff{n-1}$.
Because
restriction and taking generalized eigenspaces are exact functors,
we have the following claim.


\begin{claim}
\label{claim-oplus}
$\ei$ is an exact functor $\Rep_q H_n \to \Rep_q H_{n-1}$
and $ \Res_{H_{n-1}}^{H_n} M = \bigoplus_{i \in \Z} \ei M$.
\end{claim}

Define
$$\epsi(M) = \max \{ m \ge 0 \mid  \ei^mM \neq 0 \}. $$
We also define
$$\eti M = \soc \ei M.$$
Observe $\epsi(M) = \max \{ m \ge 0 \mid  \eti^mM \neq 0 \} $ as well.
We remind the reader that 
the socle of a module $N$, denoted $\soc(N)$,
is the largest semisimple submodule of $N$, and that the cosocle of $N$,
denoted $\cosoc(N)$, is its largest semisimple quotient.

We recall from \cite{GV} Theorem B.
\begin{thm}[\cite{GV} Theorem B] 
\label{thm-B}
Let $M$ be an irreducible $\Haffn$-module.
Then $\eti M$ is zero or is irreducible,
and $\cosoc \ei M$ is isomorphic to $\eti M$.
\end{thm}

Similarly, given $M \in \RqHaffn$, define $\ehati M$ to
be the $X_1 - \qi$ generalized eigenspace of
$\Res_{1,n-1}^n M$, further restricted to $\Haff{n-1}$.
(Technically, one should then re-index the $T_k$ and $X_k$ to
$T_{k-1}$ and $X_{k-1}$.)
If $M$ is irreducible, let $\ethati M = \soc \ehati M$
and $\epshati(M) = \max \{ m \ge 0 \mid  \ehati^mM \neq 0 \}. $
Analogous to theorem  \ref{thm-B}, if $M$ is irreducible, then 
$\ethati M$ is irreducible or zero.
\subsection{Multisegments}
\label{sec-defnmulti}
Write  $\del ij = (\qq ij)$, $i \le j$
for the one-dimensional
 trivial representation of $\Haff{j-i+1}$ 
 on which each $T_k - q$ and $X_k - q^{k + i-1}$ vanishes.
We refer to
$\del ij$ as a {\it segment\/}.
Adopt the convention that $\del j{j-1} = \1$, the one dimensional
$\Haff{0}$-module, or that $\{ \del j{j-1} \} = \emptyset$.

We introduce the symbol $\0$, not to be confused with $\emptyset$,
which will stand for the zero module.

Call a multiset of segments a {\it multisegment\/}.
Theorem \ref{thm-class} below is the result from \cite{BZ,Z} that
multisegments classify the irreducible modules of $\Rq$.

In the following theorems, we will consider two 
 total orderings on segments. 
Examples will be given in section \ref{sec-ex}.
\begin{description}
\item[\lex order] Order segments so
$\del {i_1}{j_1} > \del {i_2}{j_2}$ if $i_1 > i_2$ or if
        $i_1 = i_2$ and $j_2 > j_1$.
\item[\std order] 
Order segments so
$\del {i_1}{j_1} \succ \del {i_2}{j_2}$ if $j_1 > j_2$ or if
        $j_1 = j_2$ and $i_2 > i_1$.
\end{description}

Here we define a function from multisegments to multisegments $\cup \{ \0 \}$
$$ \Del \mapsto \E{j} \Del.$$

\begin{rul}
\label{defn-Ei}
 Given a multisegment $\Del$, first put $\Del$ in \lex order,
so that $\Del =
\{ \del{a_1}{b_1},$   $  \del{a_2}{b_2}, \ldots,  \del{a_m}{b_m} \}$
with $a_1 \ge a_2 \ge \cdots \ge a_m$.
We write its ``$j$-signature'' as follows.
Assign a blank to any segment such that $b_k \neq j, j-1$,
assign $-$ if $b_k = j$, and   $+$ if $ b_k = j-1$.
In the corresponding word
$ \pm \pm \cdots \pm$ (interspersed with blanks), cancel any adjacent
$-+$, continuing, ignoring all previously cancelled symbols
until what is left uncanceled has the form
$$++ \cdots +-- \cdots -.$$
Suppose that the leftmost uncanceled $-$ (if it exists)
is in position $i$, i.e.~ that it came from
$\del{a_i}{j}$.
Define
$$\E{j} \Del = \Del \cup \{\del{a_i}{j-1} \} \setminus \{ \del{a_i}j \},$$
and $\E{j} \Del = \0$ if no such $-$ exists. 
  \end{rul}

Observe that $\pm$ word corresponding to $\E{j} \Del$,
if $\E{j} \Del \neq \0$, is exactly that of
$\Del$ except that the leftmost uncanceled $-$ has been changed to $+$.
Further, if $\ve$ is the total number of uncanceled $-$ signs,
then $\E{j}^{\ve +1} \Del = \0$.


Define a second function from multisegments to multisegments $\cup \{ \0 \}$
$$ \Del \mapsto \Ehati \Del.$$

\begin{rul}
\label{defn-Ehat}
Given a multisegment $\Del$, we define $\Ehati \Del$ as follows:
put $\Del$
into \std order.
For each $\del iz$ write  $-$, for each $\del{i+1}{z'}$ write  $+$,
and for all other segments write a blank.
Now we ignore all $+-$ pairs, similar to rule \ref{defn-Ei},
which will leave uncanceled symbols
$- \cdots -+ \cdots +$.
Then, if $\del{i}z$ corresponds to the rightmost $-$, we replace that segment
with $\del{i+1}z$.
\omitt{In other words, the rightmost $-$ gets changed to a $+$,
and $\Ehati \Del = \Del \setminus \{\del{i}z\} \cup \{\del{i+1}z\}$.}
If no such uncanceled $-$ exists, then $\Ehati \Del =\0$.
\end{rul}

   Given a multisegment and a choice of ordering of its segments
   $\Del = \{\del{i_1}{j_1}, \ldots,$  $ \del{i_m}{j_m} \}$ let
   $$\Ind \Del$$
   denote $\Ind_{n_1, n_2, \ldots, n_m}^n \del{i_1}{j_1} \bx \cdots
   \bx \del{i_m}{j_m}$,
   where $n_k = j_k - i_k +1$ and $n = \sum_k n_k$.
   Define $n(\Del) = n$ and $m(\Del) = m$.
(Observe $n( \E{j} \Del) = n(\Del) -1$ if $\E{j} \Del \neq \0$.)

   Multisegments are by definition unordered.  However, it 
   is necessary when discussing induced modules to have an order understood.
   From now on, unless stated otherwise,  given a multisegment $\Del$, let
   $\xDel$ denote the multisegment in \lex order, and $\sDel$
   denote it in \std order.

Define $$\Md = \cosoc \Ind \Del.$$
Set  $\Mg{\0} = 0$.
Again, the symbol $\0$ is not to be confused with the empty segment $\emptyset$
for which $\Mg{\emptyset} = \1$, the one dimensional module
of $\Haff{0} = R$.

The following theorem of \cite{BZ, Z} show that multisegments
parameterize the irreducible modules in $\RqHaff$, $n \ge 0$.
A proof will be given in section \ref{sec-class}.
   \begin{thm}[\cite{BZ,Z}]
   \label{thm-class}
   Let $\Del$ be a multisegment with $n = n(\Del)$. Then
    \begin{enumerate}
   \item
   $\Md  := \cosoc \Ind \xDel$ is an irreducible $\Haffn$-module.
   \item
   If $\Del \neq \Del'$ then $\Md \not\isom \Mg{\Del'}$.
   \item
   Given any irreducible $M \in \RqHaffn$, there exists a multisegment
   $\Del$ (with $n(\Del) = n$) such that $M = \Md$.
   \end{enumerate}
   \end{thm}

The main theorem of section \ref{sec-main}, proved in section
\ref{sec-proofseg}, is that
$\et{j} \Md = \Mg{\E{j} \Del}$
(and   $\ethati \Md = \Mg{\Ehati \Del}$).

\subsection{Partitions}
\label{sec-defnpart}
Recall that $\mu$ is a partition of $n$ with length $k$ if 
$\mu = (\mu_1, \ldots, \mu_k)$ such that $\sum_m \mu_m = n$
and $\mu_1 \ge \cdots \ge \mu_k >0$. We write $|\mu| = n$.

\begin{defn}
\label{defn-Del-of-mu}
For a partition $\mu$ of length $k$ and for  $i \in \Z$, let
$\Del(\mu, i) = \{ \del{i}{i+\mu_1-1}, \ldots,$
$ \del{i - k+1}{i-k +\mu_{k}} \}$.
Observe $|\mu|  = n(\Del)$.
\end{defn}

We will say $\mu$ is colored by $i$ if we associate it
with $\Del(\mu, i)$. 
We picture it as follows.
Consider the Young diagram associated
to $\mu$,
which consists of $\mu_1$ boxes in the first row,
$\mu_2$ boxes in the second, etc.
Fill the $(x,y)$ box with $i+x-y$.
$$ \tableau{i& {\scriptstyle i+1} 
& \cdots&\mbox{}& \cdots& b_1 \\
      {\scriptstyle i-1} & i& \cdots& b_2\\
      \vdots& \mbox{}&\mbox{}\\
      {\scriptscriptstyle i\!-\!k\!+\!1} & \cdots & b_{k}}$$
Then the row fillings correspond to the segments of $\Del(\mu, i)$),
with $b_m = i + \mu_m -m$.
Then we say $\mu,i$ has a {\it removable\/} $j$-box if
you can remove a $j$-filled box $\jboxj$ from the diagram such that the
result is again the diagram of a partition (colored by $i$).
We say $\mu,i$ has an {\it addable\/} $j$-box if
you can add a $j$-filled box to the diagram such that the
result is again the diagram of a partition (colored by $i$).

We will call an $r$-tuple of partitions 
$\underline{\mu} = (\mu^{(1)}, \ldots, \mu^{(r)})$ a multipartition of $n$
if $n = \sum_m |\mu^{(m)}|$.

Likewise, given a multipartition $\underline{\mu}$ and a weight
$\lambda = \Lambda_{i_1} + \cdots + \Lambda_{i_r}$ of level $r$,
with $i_1 \le i_2 \le \cdots \le i_r$,
we will say $\underline{\mu}$ is 
colored by $\lambda$ if  $\mu^{(m)}$ is colored by $i_m$.

%

\label{defn-Kmp}
If $\mubarlam$ satisfy the condition
\begin{gather*}
\mu_{i_t - i_{t+1} + x}^{(t)} \le \mu_x^{(t+1)} \quad
\text{ for all $x \ge 1$, $1 \le t \le r-1$}
\end{gather*}
then we will say it is a Kleshchev multipartition.
The only multipartitions we will consider in Theorems \ref{thm-part1}
and \ref{thm-part2} will be Kleshchev multipartitions.
See corollary \ref{cor-part4tensor} for a crystal-theoretic interpretation
of Kleshchev multipartitions.

\begin{defn}
\label{defn-Del-of-mubar}
If $\mubarik$ is colored by $\lambda$, set
$\Del(\mubarlam) = \bigcup_{j=1}^r \Del(\mu^{(j)},i_j)$.
\end{defn}

\begin{rul}
\label{defn-p}
Given a 
multipartition $\underline{\mu}$
and weight $\lambda = \Lambda_{i_1} + \cdots + \Lambda_{i_r}$
with $i_1 \le i_2 \le \cdots \le i_r$
 we define $\E{j} (\mubarlam)$
as follows.  If $\mu^{(k)}$
 has a removable $j$ box when colored by $i_k$,
write the symbol $-$.
If  $\mu^{(k)}$ has an addable $j$ box when colored by $i_k$, write $+$. 
Otherwise,
write a blank.  In the resulting word of length $r$,  cancel
all occurrences of $-+$.
In the remaining uncanceled symbols
$++ \cdots +- \cdots -$, we remove the $j$ box from the  $\mu^{(k)}$
that corresponds to the leftmost $-$ symbol, if one exists.
The resulting colored multipartition is $\E{j} (\mubari)$.
Otherwise, $\E{j} (\mubari) = \0.$
\end{rul}

Again, we note that $\E{j}$ ``changes'' the leftmost $-$ to a $+$.

\begin{defn}
\label{defn-Nm}
Define $\Nmui = \Mg{\Del(\mu, i)}$. 
Let the symbol $\0$ be such that 
$\Del(\0, i) = \0$ and
$\Nm \0i = 0$ (not to be confused with the empty partition
$\mu  = \emptyset$, in which case $\Nm{\emptyset}{i} = \1$).
Define $\Nmulam =
\cosoc \Ind \Nm{\mu^{(1)}}{i_1} \bx \cdots \bx \Nm{\mu^{(r)}}{i_r}$.
\end{defn}

In Theorem \ref{thm-part1}, we will see the
modules $\Nmulam$ are irreducible when $\mubarlam$ is a 
Kleshchev multipartition.
In section \ref{sec-main}, we will also see the rule for $\E{j}$ on colored
multipartition is compatible with that for the corresponding multisegment and
with that of $\et{j}$ on the corresponding irreducible module.

\section{Main results}
\label{sec-main}

Below are the theorems that  give the action of $\et{j}$ on irreducible
modules of $\Rq$ both in terms of parameterization by
multisegments and by colored multipartitions.
  The action of $\et{j}$ on irreducibles corresponds to the
action of $\E{j}$ on multisegments and on Kleshchev 
multipartitions.
In other words, of the three crystal graphs 
$$\B_{BZ}, \qquad \qquad \B_\aff, \qquad \qquad \B'_\lambda$$
with edges described by 
$$
\Ei \Del  \xrightarrow{i} \Del
\qquad
\eti M \xrightarrow{i} M
  \qquad
\Ei (\mubarlam) \xrightarrow{i} \mubarlam
$$
the first two are isomorphic with isomorphism given by 
$$\Del \mapsto \Md.$$
The third is a connected component of the graph of {\it all\/}
$\lambda$-colored multipartitions, and is isomorphic
to the subgraph $\B_\lambda$
of $\B_\aff$ corresponding to $\Rep \Hlamn$ for $n \ge 0$.
That embedding is
$$
\mubarlam \mapsto \Nmulam \quad \text{ or } \quad
\mubarlam \mapsto \Del(\mubarlam).$$

Proofs are postponed until sections
\ref{sec-proofseg} and \ref{sec-proofpartition}.
   \begin{thm}
\label{thm-Ei}
Let $\Del$ be any multisegment. Then $\et{j} \Md = \Mg{\E{j} \Del}$.
   \end{thm}

We have the analogous theorem
   \begin{thm}
\label{thm-Ehat}
Let $\Del$ be any multisegment.
Then $\ethati \Md = \Mg{\Ehat{i} \Del}$.
   \end{thm}

As a corollary, we determine which $\Md$ are in $\RHlamn$.
\begin{cor}
\label{cor-cycl}
Let $\Del$ be a multisegment with $n = n(\Del)$.
Let $\lambda = \sum m_i \Lambda_i$.
Then $\Md \in \RHlamn$ $\iff$ for all $i$,
in the $\pm$ word computed as in Rule \ref{defn-Ehat}
to calculate
$\Ehati \Del$,
there are $\le m_i$ uncanceled $-$ signs.
In other words, $\Ehati^{m_i +1} \Del = \0$.
\end{cor}

\begin{proof}
This follows from theorem \ref{thm-Ehat} and also 
 relies on Theorem 9.13 of \cite{G}: that 
$\epshati(M)$ is the maximal size of the Jordan block
with eigenvalue $\qi$ for $X_1$ acting on an irreducible
module $M$.
\end{proof}

Theorem \ref{thm-Ehat} describes a different isomorphism
$\B_{BZ} \to \B_\aff$ than Theorem \ref{thm-Ei},
and hence an automorphism of $\B_\aff$.
This automorphism is also given by 
$$M \mapsto \rev^* M$$
where $\rev:\Haffn \to \Haffn$ is the algebra automorphism
$\rev(X_k) = X_{n+1-k}$, $\rev(T_k) = -(T_{n-k} + 1 - q)$.
In fact, it is this involution we really use to 
deduce Theorem \ref{thm-Ehat} from Theorem \ref{thm-Ei}.
Because $\rev$ exchanges $X_1$ and $X_n$, any statement 
made regarding $\ei$ and an irreducible module $M$
can be made for $\ehati$ and $\rev^* M$.


\begin{thm}
\label{thm-part1}
Let $\Md$ be an $\Hlamn$-module where $\lambda = \Lambda_{i_1} +
\Lambda_{i_2} + \cdots + \Lambda_{i_r}$, $i_1 \ge \cdots \ge i_r$.
Then there exists an $r$-tuple of partitions
$\underline{\mu} = (\mu^{(1)}, \ldots, \mu^{(r)})$ satisfying
\begin{gather}
\label{eq-connected1}
\mu_{i_t - i_{t+1} + x}^{(t)} \le \mu_x^{(t+1)} \quad
\text{ for all $x \ge 1$, $1 \le t \le r-1$}
\end{gather}
such that $\Md =
\Nmulam$
and $\Del =
\Del(\mubarlam)$.

Conversely, if $\mubari$ satisfies \eqref{eq-connected1}, then
$\Nmulam$
$=\cosoc \Ind \Nm{\mu^{(1)}}{i_1} \bx \cdots \bx \Nm{\mu^{(r)}}{i_r}$
is irreducible and isomorphic to $\Md$ for $\Del=
\Del(\mubarlam)$.

 \end{thm}

Recall that we call an $r$-tuple of partitions satisfying
 condition \eqref{eq-connected1} a Kleshchev multipartition.
\begin{thm}
\label{thm-part2}
For a colored multipartition satisfying \eqref{eq-connected1},
$\et{j} \Nmb{\mubarlam} = \Nmb{ \E{j} (\mubarlam)}$.

Further, $\E{j}( \Del(\mubari) ) = \Del( \E{j} (\mubarlam))$.
\end{thm}

\begin{cor}
\label{cor-part3}
For fixed $\lambda$, and for $\mubarik$ satisfying \eqref{eq-connected1},
the modules
$$\Nmulam = \cosoc \Ind \Nm{\mu^{(1)}}{i_1} \bx \cdots \bx \Nm{\mu^{(r)}}{i_r}$$
are irreducible, distinct, and range over all irreducible
$\Hlamn$-modules.
In other words, Kleshchev multipartitions parameterize the
irreducible $\Hlamn$-modules.
\end{cor}

It is instructive to compare this result and the following result
with similar results
of \cite{AM, A}.
Both their construction of  the modules $\Nmulam$ and the proof
that they are nonzero are different from ours.

In section \ref{sec-crystals},
we explain the following corollary, which  interprets the preceding
theorems in terms of the crystals $\B(\lambda)$, $\B(\infty)$ and
$\B(\Lambda_{i_r}) \otimes \cdots \otimes \B(\Lambda_{i_1})$.

\begin{cor}
\label{cor-part4tensor}
The crystal graph $\B_\lambda'$ is isomorphic to $\B_\lambda$
via 
$$\mubarlam \mapsto \Nmulam$$
and is isomorphic to
the connected component of the unique node of weight $\lambda$ in
$\B(\Lambda_{i_r})  \otimes \cdots \otimes \B(\Lambda_{i_1})$ 
via
$$ \mubarlam \mapsto \mu^{(r)} \otimes \cdots \otimes \mu^{(1)}.$$
\end{cor}

\section{Further definitions}
\label{sec-def2}

We need a few more definitions to proceed. Here they are.

If $A$ is an $R$-algebra, we write
$\Rep A$ for the category of left $A$-modules which are
{\it finite dimensional\/} as $R$-modules.
We recall that the socle of a module $M$, denoted $\soc(M)$,
is the largest semisimple submodule of $M$, and that the cosocle of $M$,
denoted $\cosoc(M)$, is its largest semisimple quotient.
We also write $Z(A)$ for the center of the algebra $A$.

If $A$ and $A^\prime$ are two $R$-algebras
%
with modules $M$ and $M'$ respectively, 
let $M \boxtimes M'$ denote the
$A \otimes A'$-module which is isomorphic to
$M \otimes_R M'$ as an $R$-module and has $A \otimes A'$
action given by
$$
(a \otimes a')\cdot (m \otimes m') = am \otimes a'm'.
$$
Because $R$ is algebraically closed, if
$M$ is an  irreducible $A$-module
and $M^\prime$ is an  irreducible $A^\prime$-module,
then 
$M \boxtimes M^\prime$ will be an
 irreducible representation of
$A \otimes A^\prime$, and all such are of this form.

\subsection{Induction and Restriction}
\label{sec-indres}

Recall that if $A \subset B$ are $R$-algebras
such that $B$ is finitely generated both as a left
and right $A$-module, the exact functor
of restriction
$$\Res_A^B :  \Rep {B} \to  \Rep {A} $$
has left and right adjoints, $\Ind$ and $\Indhat$ defined by
$$\Ind_A^B :  \Rep {A} \to \Rep {B}   \hspace{3ex}
M \mapsto B \otimes_A M$$
$${\Indhat}_A^B :  \Rep {A} \to \Rep {B}   \hspace{3ex}
 M \mapsto \Hom_A(B, M);$$
i.e.
$$\Hom_B(\Ind M, N) = \Hom_A (M, \Res N) \quad
\Hom_B(N, {\Indhat}M) = \Hom_A(\Res N, M).$$
If $B$ is a free $A$-module, then $\Ind$ and $\Indhat$
are exact functors also.
Further, if $A \subset B \subset C$ are inclusions of $R$-algebras,
we have transitivity of induction and restriction:
$$\Res_A^B \Res_B^C = \Res_A^C, \hspace{3ex}
\Ind_B^C \Ind_A^B = \Ind_A^C,    \hspace{3ex}
\Indhat_B^C  \Indhat_A^B = \Indhat_A^C$$

Now apply these remarks to the affine Hecke algebra.
Given a sequence $P= (\aaa_1, \ldots, \aaa_k)$ of non-negative
integers summing to $n$, write
$H_P = H_{\aaa_1} \otimes \cdots \otimes H_{\aaa_k}$.
We have an obvious embedding
$$H_P = H_{\aaa_1} \otimes \cdots \otimes H_{\aaa_k}
\hookrightarrow H_n$$
which makes
$H_n$
a free $H_P$-module.
Applying the previous remarks we get exact functors
$\Res, \Ind, {\Indhat}$.
When unambiguous from context, we just write
$\Res_P^n$ or $\Res$ for $\Res_{H_P}^{H_n}$, and similarly
for $\Ind$ to lighten notation.
These functors depend on the order $(\aaa_1, \ldots, \aaa_k)$ and not
just on the underlying set. 

%
%
%
%

\subsection{characters}
\label{sec-ch}

Write $S= \Rx$ and let $M$ $\in \RqHaffn$.

Define the generalized $S$-eigenspace
 $$M[\gamma] = \{v \in M \mid (X_i - \gamma_i)^m v = 0, \,
m \gg 0, \text{ for all } 1 \le i \le n \}$$
where $\gamma = (\gamma_1, \ldots, \gamma_n) \in R^n.$

Because the $X_i$ commute, 
we can decompose $M$ into a direct sum of generalized
eigenspaces
$$\Res_S^{\Haffn} M = \bigoplus_{\gamma \in \{\qi\}^n} M [\gamma].$$

We define the character of $M$ to be
the formal sum
$$\ch M = \sum_{\gamma \in (R^\times)^n} (\dim M[\gamma]) \gamma.$$

Since $S =\Rx = R[X_1^{\pm1}] \otimes_R \cdots \otimes_R R[X_n^{\pm1}]
= \Haff{(1,1,\ldots,1)}$,
we will  also write
$ \gamma = \gamma_1 \boxtimes \gamma_2 \boxtimes \cdots \boxtimes \gamma_n$
for the
1-dimensional $S$-module with character $\gamma = (\gamma_1, \ldots, \gamma_n)$.
For ease, we will also write $\gamma = \gamma_1 \gamma_2 \cdots \gamma_n$
for a term in the character of an $\Haffn$-module.

\section{Useful propositions}
\label{sec-props}
In this section, we collect useful results that we will need later.
The results in this section are well-known and most can
be found in \cite{B}, \cite{G}, \cite{GV}, \cite{Ka}, \cite{V2}. 
 
\begin{prop}
If $M$ is an irreducible $\Haffn$-module, then $M$ is finite
dimensional  with $\dim_R M \le n !$.
\label{prop-dim}
\end{prop}
Consequently, $\Rep \Haffn$ includes all irreducible $\Haffn$-modules.
\begin{thm}[Kato's Theorem]
\label{thm-kato}
Let $ \qi K_n$ denote   $\Ind_{1,1,\ldots,1}^{n} 
\qi \boxtimes \cdots \boxtimes  \qi$.
Then $\qi K_n$
is irreducible, and  is the unique $H_n$-irreducible module
with
$\qi \qi \cdots \qi$ occurring in its character.
\end{thm}
From \cite{G} Proposition 12.1 we have:
\begin{prop}[Serre relations]
 \label{prop-serre}
As operators on the Grothendieck group $\oK$
\begin{enumerate}
\item if $| i-j| > 1,$ then $\ei \e{j} = \e{j} \ei$.
\item $\ei^2 \e{i \pm 1} + \e{i \pm 1} \ei^2
= 2 \ei \e{i \pm 1} \ei$.
\end{enumerate}
\end{prop}
%
\begin{rem}
\label{rem-serre}
We can reinterpret proposition \ref{prop-serre},
that, first, 
if $q^j \neq q^{i \pm 1}$, then for   any finite dimensional module $M$,
if $\gamma_1 \gamma_2 \cdots q^i q^j \cdots \gamma_n$ occurs in $\ch M$
then so does 
  $\gamma_1 \gamma_2 \cdots q^j q^i \ldots \gamma_n$.
(Specifically:
$\dim M[(\gamma_1, \gamma_2, \ldots, q^i, q^j, \ldots, \gamma_n)] = $
  $\dim M[(\gamma_1, \gamma_2, \ldots, q^j, q^i, \ldots, \gamma_n)]$.)

Secondly, 
$\gamma_1  \cdots \qi q^{i \pm 1} \qi \cdots \gamma_n$ occurs in $\ch M$
if and only if
$\gamma_1 \cdots \qi \qi  q^{i \pm 1} \cdots \gamma_n$
{\it or\/}
$\gamma_1 \cdots   q^{i \pm 1} \qi \qi \cdots \gamma_n$
occur.
(Furthermore,
the dimensions of the  generalized $S$-eigenspaces spaces satisfy that
twice the first is the sum of the next two.)

Kato's Theorem implies that if 
$\gamma_1  \cdots \underbrace{\qi  \qi  \cdots \qi}_k \cdots \gamma_n$
occurs in $\ch M$, it does so with multiplicity at least $k!$.
\end{rem}

\begin{prop}[\cite{GV},  Lemma 3.5, part 3]
\label{lemma-3.5}
For any irreducible module $N$
$$M = \cosoc \Ind N \bx \qi K_k$$
is irreducible with $\epsi(M) = \epsi(N) + k$, and
all other composition factors of $\Ind N \bx \qi K_k$
have strictly smaller $\epsi$.
\end{prop}

\begin{prop}[\cite{GV},  Corollary 3.6]
\label{cor-3.6}
Let $M$ be  an irreducible $H_n$-module.
\begin{enumerate}
\item
Let $\ve = \epsi(M)$ and $N = \eti^\ve M$.
Then for $0 \le k \le \ve$
\begin{gather*}
\eti^k M = \cosoc \Ind N \boxtimes  \qi K_{\ve -k}.
\end{gather*}
\item \label{chain-2}
$\eti M \isom \cosoc \ei M$.
\item If $L$ is also irreducible and $\eti M \isom \eti L \neq 0$
then $M \isom L$.
\end{enumerate}
\end{prop}

\begin{cor}
\label{cor-e-et}  
Let $M$ be an irreducible $\Haffn$-module and let  $\ve = \epsi(M)$.
Then $\ei^\ve M$ is the direct sum of $\ve !$ copies of
$\eti^\ve M$.
\end{cor}
\begin{proof}
From proposition \ref{cor-3.6}, if $N = \eti^\ve M$, then
$M = \cosoc \Ind_{n-\ve, 1, \ldots, 1}^n N \boxtimes q^i 
\boxtimes \cdots \boxtimes q^i=
 \cosoc \Ind_{n-\ve, \ve}^n N \boxtimes q^i K_\ve$
and $\epsi(N)=0$.
Applying the exact functor $\ei^\ve$  (and using lemma \ref{lemma-mackey})
yields
a surjection 
\begin{gather}
\label{eq-ve}
\bigoplus_1^\ve N \surj \ei^\ve M.
\end{gather}
But by Frobenius reciprocity, we have a map
$N \bx \qi K_\ve \to \Res_{n-\ve, \ve}^n M$, which must
   be an injection by Kato's theorem. Since restriction is exact,
this shows
\eqref{eq-ve} must be an isomorphism.
\end{proof}

The next proposition follows from proposition
\ref{cor-3.6} and  \cite{G} Proposition 10.4.
\begin{prop} 
\label{prop-10.4}
Let $M$ be an irreducible $\Haffn$-module such that   $\epsi(M) > 0$.
Then
\begin{enumerate}
\item  $\epsi(\eti M) = \epsi(M) -1$.
\item $\epsj{i \pm 1}(\eti M) = 
\begin{cases} \epsj{i \pm 1}(M) 
	\qquad \text{ or} \\
	\epsj{i\pm 1}(M) + 1 \end{cases}$
\end{enumerate}
\end{prop}

 \subsection{The shuffle lemma}
\label{sec-shuffle}

We write $S_n$ for the symmetric group on $n$ letters,
$\ttau_i = (i \hspace{2ex} i+1)$, $1 \le i \le n-1$,
 for the simple transpositions.
Denote  length by  $\ell(w)$. 

If $P = (\aaa_1, \ldots, \aaa_k)$ is an {\it ordered \/}
tuple, $\aaa_i \in \Z_{>0}, \sum \aaa_i = n$,
it is convenient to denote $S_P = S_{\aaa_1} \times \cdots \times S_{\aaa_k}
\subseteq S_n$.
For example, $S_{(n)} = S_n$, $S_{(1,\ldots,1)} = \{ 1 \}$.

For such $P$, we write
$\WP{P}$ for the set of minimal length left coset
representatives of $S_P \subseteq S_n$.

The following lemmas are  special cases of the ``Mackey
formula'' relating the composite of induction from $H_P$
to $H_n$ with the 
restriction from $H_n$ to $H_{P'}$, for various $P$ and $P'$.
In particular, in lemma \ref{lemma-mackey} $P= (a,b)$, $P' = (n-1,1)$
(but then we further restrict to $\Haff{n-1}$);
in lemma \ref{lem-chind} below,
$P = P' = (1, 1, \ldots, 1)$ so that $H_P = S$;
and finally in lemma \ref{lemma-shuffle},
$P = (m,n)$ and $P'= (1, 1, \ldots, 1)$.

Below we compute the action of $\ei$ on an induced module. We refer the reader
to \cite{G} or \cite{GV} for a proof.
 \begin{lemma} 
 \label{lemma-mackey}
Let $A$ be an irreducible $\Haff{a}$-module and $B$ be an
irreducible $\Haff{b}$-module.  Let $n=a+b$.
The following is an exact sequence:
\begin{gather*}
0 \to  \Ind_{a-1,b}^{n-1} \ei A \bx B
 \to \ei (\Ind_{a,b}^n A \bx B) \to \Ind_{a,b-1}^{n-1} A \bx \ei B \to 0
\end{gather*}
 \end{lemma}

Theorems \ref{thm-Ei} and \ref{thm-part2} describe that for very
special $A$ and $B$, we can determine exactly when
\begin{gather*}
\eti (\cosoc \Ind A \bx B) = \begin{cases}
	\cosoc (\Ind \eti A \bx B)  \\
	\cosoc (\Ind A \bx \eti B) 
\end{cases}
\end{gather*}
and that the cosocles above are irreducible.
The first case happens when 
$\epsi(A) > \varphi_i(B)$, and the second when
$\epsi(A) \le \varphi_i(B)$.
We will define $\varphi_i$ in section \ref{sec-crystals}.

Compare this to equation \eqref{eqn-crystal2} in section \ref{sec-crystals}
concerning a tensor product of crystals.

\begin{lemma}
\label{lem-chind}
Let $\gamma = (\gamma_1, \ldots, \gamma_n) \in \{ \qi \}^n$. 
Then
$$\ch(\Ind_S^{H_n} \gamma_1 \boxtimes \cdots \boxtimes \gamma_n)
= \sum_{w \in S_n} w\cdot \gamma.$$
\end{lemma}

We now describe the character of $\Ind_P^n M$ in terms
of $\ch M$.
We say a string $\gamma = \gamma_1 \gamma_2 \cdots \gamma_k$
is a shuffle of $t$ and
$u$ if $t$ is a subword of $\gamma$ and $u$ is its complementary subword.
The shuffle of $t$ and $u$, denoted $t \shuffle u$,
is then the formal sum of all shuffles of $t$ and $u$, with multiplicity.

The permutations of $\{1, 2,  \ldots,  n+m\}$ that keep
$1,  \cdots,  m$ in order and $m+1,  \cdots,  m+n$ in order
(i.e.~their shuffles)
are exactly the minimal length left coset representatives
$\WP{(m,n)}$.
It follows that if we write $tu$ for the concatenation of $t$ and $u$, which
have length $m$ and $n$ respectively, then
$$t \shuffle u = \sum_{w \in \WP{(m,n)}} w\cdot (tu).$$


We  extend $\shuffle$ linearly to sums of words.
\begin{lemma}
\label{lemma-shuffle}
Let $M \in \Rep H_m$, $N \in \Rep H_n$. Then
$$
\ch \Ind_{(m,n)}^{{m+n}} M \boxtimes N
= \sum_{\substack{\gamma \text{ is a shuffle}\\ \text{ of $t$ and $u$}}}
(\dim M[t] \dim N[u]) \gamma = \ch M \shuffle \ch N.$$
\end{lemma}

\subsection{Induced modules}
\label{sec-ind}
\begin{prop}
\label{prop-rev}
Let $A$ be an irreducible $\Haff{a}$-module and $B$ be an
irreducible $\Haff{b}$-module.  Let $n=a+b$.
Then $\Ind_{a,b}^n A \bx B \isom \Indhat_{b,a}^n B \bx A$, where $n = a+b$.
\end{prop}
\begin{proof}
A full proof can be found in Proposition 3.5 of \cite{V2}.
We will outline the
construction of an explicit isomorphism.

Let $W = \WP{(a,b)} = \{ w \} = $ minimal length left coset
representatives for $S_a \times S_b \subseteq S_n$.
Let $W' = \{ x_0 w^{-1} \}_{w\in W}= \{ w' \} =$
minimal length right coset representatives
for $S_b \times S_a \subseteq S_n$, where $x_0$ is the
longest element of $W$.
If $\{ u \tensor v \}$ is a basis of $A \bx B$, then
$\{  v \tensor u \}$ is a basis of $B \bx A$, so
$\{ T_w \otimes (u \tensor v) \}_{w \in W}$
is a basis of $\Ind_{a,b}^{n} A \bx B$, and
$\{ \phi_{w', u \otimes v}\}_{w' \in W'}$
is a basis of $\Indhat_{b,a}^{n}{B \bx A} =
\Hom_{\Haff{(b,a)}}(\Haffn, B \bx A)$,
where $\phi_{w', u \tensor v}(T_{x'}) = \delta_{w', x'}  v \tensor u$
for $x' \in W'$ (and extend to a homomorphism on all of $\Haff{(b,a)}$).

Write $\phi_{u\tensor v}$ for $\phi_{x_0, u \tensor v}$.
Define the isomorphism
\begin{align*}
\psi : \Ind_{a,b}^{n} A \bx B \longrightarrow &\,
\Indhat_{b,a}^{n} B \bx A
%
\\
h \otimes (u \tensor v) \stackrel{\psi}{\longmapsto} & \, h \phi_{u\tensor v}
 \qquad \text{ for $h \in H_n$.}
\end{align*}
To check that $\psi$ is well-defined, it suffices by Frobenius
reciprocity to check that the map $u\tensor v \mapsto \phi_{u \tensor v}$ is
an $\Haff{a} \tensor \Haff{b}$-map.
Then one must show $\psi$ is surjective, yielding it is an isomorphism
by a dimension count.
We leave these details to the reader.
\end{proof}
\begin{prop}
\label{prop-soccosoc}
Let $N = \cosoc \Ind_{a,b}^n A \bx B$ where $n = a+b$.
Then $N \isom \soc \Ind_{b,a}^n B \bx A \isom
\soc \Indhat_{a,b}^n A \bx B$.
\end{prop}
\begin{proof}
It is well-known  that there exists an involution
$D : \Rq \to \Rq$ that takes irreducibles to themselves
and
$$0 \to A \to B \to C \to 0 \iff 
0 \to D(C) \to D(B) \to D(A) \to 0.$$
This is an unpublished result of Bernstein \cite{B} that has
since entered the literature.
This functor also has the property that
$$D(\Ind M) = \Indhat D(M).$$
In other words, there is an even stronger relationship
between $\Ind$ and $\Indhat$:  for $M$ irreducible,
the factors in the socle series of $\Ind M$
coincide with the factors in the cosocle series of $\Indhat M$
{\it in order\/}.
\end{proof}

We remark that where we use proposition \ref{prop-soccosoc}
in this paper, it is possible to use alternate arguments using characters,
but they are  far less elegant.

\begin{prop}
 \label{prop-ch}
 Let $P = (n_1, \ldots, n_k)$ and $n = \sum_i n_i$.
 Let $\gamma_i$ be a one-dimensional $H_{n_i}$-module,
 and so we can write $\gamma_i$ for its character as well.
(In section \ref{sec-ch} we took $n_i = 1$.)
If $Q$ is any quotient of $\Ind_P^n  
\gamma_1 \bx \cdots \bx \gamma_k$,
then $\ch Q$ contains the concatenation $ \gamma_1 \gamma_2 \cdots \gamma_k$.
If $L$ is any submodule of $\Ind_P^n  
\gamma_1 \bx \cdots \bx \gamma_k$
then $\ch L$ contains the term
$\gamma_k \cdots \gamma_2 \gamma_1$.          
  \end{prop}
\begin{proof}
Frobenius reciprocity gives $\Hom_{\Haffn}(\Ind_P^n 
\gamma_1 \bx \cdots \bx \gamma_k, Q)
= \Hom_{\Haff{P}}(\gamma_1 \bx \cdots \bx \gamma_k, \Res_P^n Q)$,
from which the first statement is immediate once we restrict
from $\Haff{P}$ to $S = \Haff{(1,\ldots, 1)}$.

Similarly, the second follows by proposition \ref{prop-rev}
since 
$\Hom_{\Haffn}(L, \Ind_P^n \gamma_1 \bx \cdots \bx \gamma_k)
= \Hom_{\Haffn}(L, \Indhat_{P'}^n \gamma_k \bx \cdots
	\bx \gamma_2 \bx \gamma_1)
= \Hom_{\Haff{P'}}(\Res_{P'}^n L, \gamma_k \bx \cdots \bx
	\gamma_2 \bx \gamma_1).$
\end{proof} 
\begin{prop}
 \label{prop-irred}
Let $M$ be an irreducible $H_P$-module.  Suppose $\gamma$ occurs as a term in
$\ch (\Ind_P^n 
M)$ with multiplicity $m$ and
also that $\gamma$ occurs  with multiplicity $m$ in the character
of any quotient of $\Ind_P^n  
M$.
Then $\cosoc \Ind_P^n 
M$
is irreducible and occurs with multiplicity one as a composition
factor of $\Ind_P^n  
M$.
If in addition  $\gamma$ occurs in $\ch (\soc \Ind_P^n   
M)$
then $\Ind_P^n  
M$ is irreducible.
\end{prop}
\begin{proof}
The first statement follows since the map 
$$\Ind M \surj \cosoc \Ind M$$
is also a homomorphism of $S = \Rx$-modules.
Restricting to $S$, the $\gamma$-isotypic component (generalized eigenspace)
of $\Ind M$ must surject to that of $\cosoc \Ind M$ because  
taking generalized eigenspaces is exact.
Hence $\gamma$ can occur with multiplicity at most $m$ in
the cosocle.
By hypothesis, each irreducible component of $\cosoc \Ind M$ accounts
for $m$ copies of $\gamma$, and so the cosocle must be irreducible.
A similar  counting argument shows it must occur with multiplicity one
as a composition factor of $\Ind M$.

Further, since the cosocle accounts for all $m$ copies of $\gamma$,
if $\gamma$ also occurs in $\ch (\soc \Ind    M)$,
then $\soc \Ind    M \supseteq \cosoc \Ind M$.
Hence, $\cosoc \Ind M$ is a submodule and we can consider the
quotient $\Ind M / (\cosoc \Ind M)$ whose character cannot have
any copies of $\gamma$ and so must be zero.
Thus $\Ind M = \cosoc \Ind M$ is irreducible (by the first statement).
\end{proof}

\begin{rem}
\label{rem-p9}
The same argument as above also shows that if
$\cosoc \Ind M$ is irreducible, occurs with multiplicity
one as a composition factor of $\Ind M$, and also occurs
in the socle, then $\Ind M$ is irreducible.
\end{rem}
 \section{More theorems on multisegments} \label{sec-bzseg}
This section is devoted to recalling results from \cite{BZ,Z}
on multisegments.
For convenience to the reader, we give complete proofs.
We end by giving the significance of these results in
the language of crystal graphs.

Recall if $i \le j$ then $\del ij = ( \qq ij)$ 
 is the one dimensional $\Haff{j-i+1}$-module on which all
$T_k - q$ and all $X_k - q^{k+i-1}$ vanish.

Observe that 
$\e{j} \del ij = \et{j} \del ij =  \del{i}{j-1}$,
(with the convention $ \del i{i-1} = \1$), $\ei \1 = 0$,
and if $k \neq j$ then $\e{k} \del ij = \et{k} \del ij =  0$.

 \subsection{Linking Rule} \label{sec-link}

 \begin{lemma}[\cite{BZ,Z} The linking rule]
\label{lemma-link}
 \begin{enumerate}
 \item[(i)\label{link1}]                        
 $\Ind \del{i}{ j} \boxtimes \del{k}{ l}$
 is irreducible if $j+1 < k$.

 \item[(ii) \label{link2}]
 $\Ind \del{i}{ l} \boxtimes \del{j}{ k}$
 is irreducible if $i \le j \le k \le l$.
 \item[(iii)]
 \label{link3}
 $N = \cosoc \Ind \del{j}{l} \bx \del{i}{ k}$ is irreducible
 if $i < j, k < l, j \le k+1$,
 and the following
 \begin{gather}
 \label{eqn-link}
 0 \to \Ind \del{j}{ k} \boxtimes \del{i}{ l} \to \Ind \del{j}{ l} \boxtimes \del{i}{ k}
 \to N \to 0
 \end{gather}
 is exact.        
 In this case, 
we say $\del{i}{k}$ and $\del{j}{ l}$
 are {\it linked\/} or are a {\it linked pair\/}.
 If $j= k+1$ then $\del{i}{ k}$ and $\del{j}{ l}$ are {\it adjacent\/}.
 \end{enumerate}
 \end{lemma}

 \begin{proof}
(i)
This is the case that the intervals $[ij]$ and $[kl]$ are far apart,
so that their union is not again an interval.

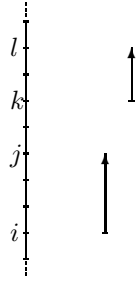
\begin{picture}(120,130)(-100,0)
\put(-30,112){\line(0,1){1}} 
\put(-30,114){\line(0,1){1}}
\put(-30,116){\line(0,1){1}}
\put(-30,14){\line(0,1){1}}
\put(-30,16){\line(0,1){1}}
\put(-30,18){\line(0,1){1}}
\put(-30,20){\line(0,1){90}} 
\put(-31,20){\line(1,0){2}}  
\put(-31,30){\line(1,0){2}}
\put(-31,40){\line(1,0){2}}
\put(-31,50){\line(1,0){2}}
\put(-31,60){\line(1,0){2}}
\put(-31,70){\line(1,0){2}}
\put(-31,80){\line(1,0){2}}
\put(-31,90){\line(1,0){2}}
\put(-31,100){\line(1,0){2}}
\put(-31,110){\line(1,0){2}} 
%
\put(0,30){\vector(0,1){30}}
\put(-1,30){\line(1,0){2}}
\put(-36,27){$i$}
\put(-36,57){$j$}
\put(10,80){\vector(0,1){20}}
\put(9,80){\line(1,0){2}}
\put(-36,77){$k$}
\put(-36,97){$l$}
\end{picture}

 Let $M = \Ind_{j-i +1, l-k+1}^n \del{i}{ j} \bx \del{k}{ l}$
 where $n = j+l-i-k+2$.
 By the shuffle lemma,
 \begin{gather}
 \label{eq-ch}
 \ch M = \sum_{w\in \WP{(j-i +1, l-k+1)}}
w \cdot (q^i \cdots q^j q^k \cdots q^l).
 \end{gather}
 Frobenius reciprocity implies that the term
 $q^i \cdots q^j q^k \cdots q^l$ must occur in any quotient of $M$,
 but by remark \ref{rem-serre} all of the
 $\binom{n}{j-i+1}$ terms in \eqref{eq-ch} must occur as well,
 yielding $M$ is irreducible.

(ii)
This is the case that one interval $[jk] \subseteq [il]$ is contained
in the other.

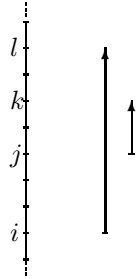
\begin{picture}(120,130)(-100,0)
\put(-30,112){\line(0,1){1}}  
\put(-30,114){\line(0,1){1}}
\put(-30,116){\line(0,1){1}}
\put(-30,14){\line(0,1){1}}
\put(-30,16){\line(0,1){1}}
\put(-30,18){\line(0,1){1}}
\put(-30,20){\line(0,1){90}}  
\put(-31,20){\line(1,0){2}}   
\put(-31,30){\line(1,0){2}}
\put(-31,40){\line(1,0){2}}
\put(-31,50){\line(1,0){2}}
\put(-31,60){\line(1,0){2}}
\put(-31,70){\line(1,0){2}}
\put(-31,80){\line(1,0){2}}
\put(-31,90){\line(1,0){2}}
\put(-31,100){\line(1,0){2}}
\put(-31,110){\line(1,0){2}}
%
\put(0,30){\vector(0,1){70}}
\put(-1,30){\line(1,0){2}}
\put(-36,27){$i$}
\put(-36,57){$j$}
\put(10,60){\vector(0,1){20}}
\put(9,60){\line(1,0){2}}
\put(-36,77){$k$}
\put(-36,97){$l$}
\end{picture}

Let $n = l-i+k-j+2$ and  $L = \Ind_{ l-i+1, k-j+1}^n \del il \bx \del jk$
where $i \le j \le k \le l$.
By  Frobenius reciprocity the term $q^i q^{i+1} \cdots q^l q^j q^{j+1}
\cdots q^k$ must occur in $\ch (\cosoc L)$.
The term $q^i \cdots q^j q^j q^{j+1} q^{j+1} \cdots  q^k q^k \cdots q^l$
must also occur in  $\ch (\cosoc L)$
by remark \ref{rem-serre} along with the shuffle lemma,
and Kato's theorem (theorem \ref{thm-kato})
implies it occurs with multiplicity at least $2^{(k-j+1)}$.
On the other hand, the shuffle lemma shows this is its multiplicity
in $\ch L$.
Similar reasoning shows this term occurs in $\ch (\soc L)$.
By proposition \ref{prop-irred}, $L = \Ind \del il \bx \del jk$
is irreducible.

(iii)
This is the case that the union of the intervals $[jl] \cup [ik]$
is a longer interval (in $\Z$). 
Let $n = l-j + k-i +2$ and  $N = \cosoc \Ind \del jl \bx \del ik$
 where $i < j, k < l, j \le k+1$.

First we'll consider the case that $j= k+1$, that is $\del jl$ and $\del ik$ are
an adjacent pair.

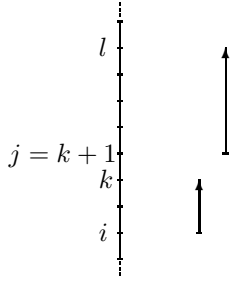
\begin{picture}(120,130)(-100,0)
\put(-30,112){\line(0,1){1}}
\put(-30,114){\line(0,1){1}}
\put(-30,116){\line(0,1){1}}
\put(-30,14){\line(0,1){1}}
\put(-30,16){\line(0,1){1}}
\put(-30,18){\line(0,1){1}}
\put(-30,20){\line(0,1){90}}
\put(-31,20){\line(1,0){2}}
\put(-31,30){\line(1,0){2}}
\put(-31,40){\line(1,0){2}}
\put(-31,50){\line(1,0){2}}
\put(-31,60){\line(1,0){2}}
\put(-31,70){\line(1,0){2}}
\put(-31,80){\line(1,0){2}}
\put(-31,90){\line(1,0){2}}
\put(-31,100){\line(1,0){2}}
\put(-31,110){\line(1,0){2}}
\put(0,30){\vector(0,1){20}}
\put(-1,30){\line(1,0){2}}
\put(-38,27){$i$}
\put(-72,57){$j =k+1$}
\put(10,60){\vector(0,1){40}}
\put(9,60){\line(1,0){2}}
\put(-38,47){$k$}
\put(-38,97){$l$}
\end{picture}
%

We want to show 
 \begin{gather}
 \label{eqn-adjacent}
 0 \to  \del{i}{ l} \to \Ind \del{k+1}{ l} \boxtimes \del i k \to N \to 0.
 \end{gather}
is exact and that $N$ is irreducible.
By Frobenius reciprocity, because $\Res_{l-k, k-i+1}^n \del i l
= \del i k \bx \del {k+1}{l}$ there is a nonzero map
$\del i l \to \Indhat_{ k-i+1,l-k}^n \del i k \bx  \del {k+1}{l} \isom$
$\Ind_{ l-k,k-i+1}^n  \del {k+1}{l} \bx \del i k$.  Further as 
$q^i \cdots q^l$
occurs only once in $\ch \Ind_{k-i+1, l-k}^n  \del {k+1}{l} \bx \del i k$,
the socle of the induced module (in \eqref{eqn-adjacent})
 must be  $\del il$.
(Note that Frobenius reciprocity (proposition \ref{prop-ch})
 shows this one-dimensional submodule
cannot be in the cosocle.)
The quotient by this submodule has character
$\sum_{\stackrel{w \in \WP{(k-i+1,l-k+1)}}{ w \neq 1 \hfill}} w \cdot \q i l$,
any one term of which (in particular $q^{k+1} \cdots q^l q^i \cdots q^k$)
ensures the occurrence of the other terms,
by remark \ref{rem-serre}.  This implies the quotient
is irreducible and hence must equal $N$.

By the same reasoning
 \begin{gather}
 \label{eqn-adjacent2}
 0 \to N \to   \Ind \del i k \bx \del{k+1}{ l} \to  \del{i}{l} \to 0
 \end{gather}
is exact.

Now consider a linked pair that is not adjacent, so that $i < j \le k < l$.

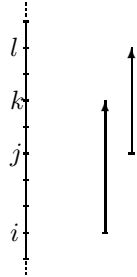
\begin{picture}(120,130)(-100,0)
\put(-30,112){\line(0,1){1}}  
\put(-30,114){\line(0,1){1}}
\put(-30,116){\line(0,1){1}}
\put(-30,14){\line(0,1){1}}
\put(-30,16){\line(0,1){1}}
\put(-30,18){\line(0,1){1}}
\put(-30,20){\line(0,1){90}}  
\put(-31,20){\line(1,0){2}}    
\put(-31,30){\line(1,0){2}}
\put(-31,40){\line(1,0){2}}
\put(-31,50){\line(1,0){2}}
\put(-31,60){\line(1,0){2}}
\put(-31,70){\line(1,0){2}}
\put(-31,80){\line(1,0){2}}
\put(-31,90){\line(1,0){2}}
\put(-31,100){\line(1,0){2}}
\put(-31,110){\line(1,0){2}}
%
\put(0,30){\vector(0,1){50}}
\put(-1,30){\line(1,0){2}}
\put(-36,27){$i$}
\put(-36,57){$j$}
\put(10,60){\vector(0,1){40}}
\put(9,60){\line(1,0){2}}
\put(-36,77){$k$}
\put(-36,97){$l$}
\end{picture}

Using \eqref{eqn-adjacent} and the exactness of induction we have an injection 
 \begin{gather}
   0 \to \Ind \del j k \bx \del i l \xrightarrow{\alpha} 
  		 \Ind \del j k \bx \del{k+1}{ l} \bx  \del i k 
\end{gather}
and from \eqref{eqn-adjacent2} a surjection
\begin{gather}
  		 \Ind \del j k \bx \del{k+1}{ l} \bx \del i k 
  \xrightarrow{\beta}   \Ind \del j l \bx \del i k  \to 0.
\end{gather}
From the shuffle lemma and the fact $i < j, k< l$, the term
$q^j \cdots q^k q^i \cdots q^j \cdots q^k \cdots q^l$ occurs
in each of the characters of
$\Ind \del j k \bx \del i l $, $\Ind \del j k \bx \del{k+1}{ l} \bx \del i k$,
and $\Ind \del j l \bx \del i k $ with multiplicity one.
Because $\alpha$ is injective and $\beta$ is surjective this implies
that $\beta \circ \alpha$ is nonzero on the $\Rx$ eigenvector with
eigenvalue $(q^j, \cdots, q^k, q^i, \cdots,  q^l)$.
By part (ii) 
$\Ind \del j k \bx \del il =
\Ind \q jk \bx (q^i \cdots q^j \cdots q^k \cdots q^l)$
is irreducible, and so $\beta \circ \alpha$ is an injection.

Observe that the term $q^j \cdots q^l q^i \cdots q^k$
occurs with
multiplicity one in $\ch(\Ind \del jl \bx \del ik)$ so that its cosocle
$N$ is irreducible by proposition \ref{prop-irred}.
 Furthermore $N$ cannot intersect $\Ind \del j k \bx \del il$ 
because this term is not in its character.
Thus we have that $\Ind \del j k \bx \del il \subseteq B$
where $B$ is the kernel of the map
\begin{gather}
  \label{eqn3} 
0 \to B \to \Ind \del jl \bx \del ik \to N \to 0.
\end{gather}
We want to show this inclusion is an equality. 

Let $A$ be any composition factor of $\Ind \del jl \bx \del ik $.
We will show $\epsj{k}(A) =1$.
First, the shuffle lemma implies $\epsj{k}(A) \le 1$.
Let $w \cdot (q^j \cdots q^l q^i \cdots q^k)$ be a term
in $\ch A$ for some $w \in \WP{(l-j+1, k-i+1)}.$
If $w(n) = n$ we are done.

If $w(n) \neq n-l+k$, then
remark \ref{rem-serre} implies the term
$s_{n-1} \cdots s_{w(n) + 1} s_{w(n)} w \cdot (q^j \cdots q^l q^i \cdots q^k)$
also occurs in $\ch A$ and this permutation does fix $n$.
In other words, the irreducible module $\Ind q^k \bx \del{l-n+w(n) +1}{l}=
\Ind \del kk \bx \del{l-n+w(n) +1}{l}
\isom \Ind \del{l-n+w(n) +1}{l}\bx \del kk$ has 
$\epsj{k} = 1$,
and this module occurs in the restriction $\Res_{n-w(n)}^{w(n), n-w(n)}
\Res_{w(n), n-w(n)}^n A$.

In the case $w(n) = n-l+k$, let $a$ be such that
 $w(n-l+k) = n-l+k +1 +a$. Then
$w \cdot (q^j \cdots q^l q^i \cdots q^k) =
 (\cdots q^k q^{k-a} \cdots q^{k-1} q^k q^{k+1} \cdots q^l)$.
Repeated application of remark \ref{rem-serre} allows us
to ``slide'' the leftmost $q^k$ over so that
$\cdots q^{k-1} q^k \cdots q^l q^k$ also occurs in $\ch A$,
yielding $\epsj{k}(A) = 1$.


Now we induct on $k-j+1$, the case $k-j+1=0$ already completed 
for adjacent pairs above. Assume $k-j+1 > 0$.
Starting with  \eqref{eqn3}, we apply the exact functor $\e{k}$ yielding
\begin{gather}
  \label{eqn4}
0 \to \e{k}B \to \Ind \del jl \bx \del i{k-1} \to \e{k} N \to 0.
\end{gather}
The middle term is $\e{k}\Ind \del jl \bx \del ik$ by
lemma \ref{lemma-mackey}.
Because $k-j+1 > 0$, the pair $\del jl, \del i{k-1}$ is still linked,
and by the inductive hypothesis we know its two composition factors.
Because $\epsj{k}(N) = 1$, 
proposition \ref{cor-3.6} implies
$\e{k} N = \et{k} N$ is irreducible, and by the inductive hypothesis
it must be equal to 
$\cosoc \Ind \del jl \bx \del i{k-1}$.
From the inductive hypothesis and the exactness of \eqref{eqn4},
we must have that $\e{k}B = \Ind \del j{k-1} \bx \del i l$,
which is irreducible by part (ii) 
of this lemma.
If $B$ were reducible, then because all composition factors of $B$ have
$\epsj{k} \neq 0$, $\e{k}B $ would also be reducible.
Hence we must have $B = \Ind \del j{k} \del i l$.
  \end{proof}

\subsection{Classification}
\label{sec-class}

Recall a multiset of segments is a {\it multisegment\/}.
Theorem \ref{thm-class} below is the result from \cite{BZ,Z} that
multisegments parameterize the irreducible modules of $\RqHaff$.

We also remind the reader of the two orderings introduced in section
\ref{sec-def1}, and include examples.

\begin{description}
\item[\lex order] 
$\del {i_1}{j_1} > \del {i_2}{j_2}$ if $i_1 > i_2$ or if
	$i_1 = i_2$ and $j_2 > j_1$.
The following picture shows the segments
$\{ \del 56, \del 57, \del 47, \del 33,$
$ \del 36, \del 36, \del 37,$
$\del 37, \del 26, \del 27, \del 29, \del {-1}7, \del {-1}1, \del{-2}2 \}$
 in \lex order,
weakly decreasing from left to right.

\begin{picture}(190,130)(-110,-10)
\put(-60,112){\line(0,1){1}}
\put(-60,114){\line(0,1){1}}
\put(-60,116){\line(0,1){1}}
\put(-60,-4){\line(0,1){1}}
\put(-60,-6){\line(0,1){1}}
\put(-60,-8){\line(0,1){1}}
\put(-60,0){\line(0,1){110}}
\multiput(-61,0)(0,10){12}{\line(1,0){2}}
\put(-75,-2){\makebox(10,5)[r]{-2}}
\put(-75,8){\makebox(10,5)[r]{-1}}
\put(-75,18){\makebox(10,5)[r]{0}}
\put(-75,28){\makebox(10,5)[r]{1}}
\put(-75,38){\makebox(10,5)[r]{2}}
\put(-75,48){\makebox(10,5)[r]{3}}
\put(-75,58){\makebox(10,5)[r]{4}}
\put(-75,68){\makebox(10,5)[r]{5}}
\put(-75,78){\makebox(10,5)[r]{6}}
\put(-75,88){\makebox(10,5)[r]{7}}
\put(-75,98){\makebox(10,5)[r]{8}}
\put(-75,108){\makebox(10,5)[r]{9}}
%
\put(-31,70){\line(1,0){2}}
\put(-30,70){\vector(0,1){10}}
\put(-21,70){\line(1,0){2}}
\put(-20,70){\vector(0,1){20}}
\put(-11,60){\line(1,0){2}}
\put(-10,60){\vector(0,1){30}}
\put(-1,50){\line(1,0){2}}
\put(0,50){\vector(0,1){4}}
\put(10,50){\vector(0,1){30}}
\put(9,50){\line(1,0){2}}
\put(20,50){\vector(0,1){30}}
\put(19,50){\line(1,0){2}}
\put(30,50){\vector(0,1){40}}  
\put(29,50){\line(1,0){2}}
\put(40,50){\vector(0,1){40}}
\put(39,50){\line(1,0){2}}
\put(50,40){\vector(0,1){40}}
\put(49,40){\line(1,0){2}}
\put(60,40){\vector(0,1){50}}  
\put(59,40){\line(1,0){2}}
\put(70,40){\vector(0,1){70}}
\put(69,40){\line(1,0){2}}
\put(80,10){\vector(0,1){80}}  
\put(79,10){\line(1,0){2}}
\put(90,10){\vector(0,1){20}}
\put(89,10){\line(1,0){2}}
\put(100,0){\vector(0,1){40}}
\put(99,0){\line(1,0){2}}
\end{picture}

\item[\std order]  
$\del {i_1}{j_1} > \del {i_2}{j_2}$ if $j_1 > j_2$ or if
	$j_1 = j_2$ and $i_2 > i_1$.
The following picture shows the {\it same \/} segments as
above in \std order,
weakly decreasing from left to right.
In other words: 
$\{ \del 29,
\del {-1}7, \del 27,$   $ \del 37, \del 37, \del 47,$ $\del 57,$
$\del 26,	\del 36,	\del 36,	\del 56,
\del 33,	\del{-2}2, \del{-1}1 \}$.
(Notice that if one rotates a picture of segments in \lex order
by $180^\circ$, one gets a picture of {\it different \/}
 segments in \std order.
However, this observation is useful in seeing why theorems for
$\ei$ also hold for $\ehati$.
In fact, if we think of the above as a picture of $\ch \Ind \Del$,
then the rotated picture (with arrows pointing down)
would be of $\ch (\rev^* M)$.)


\begin{picture}(160,130)(-80,-10)
\put(-45,-2){\makebox(10,5)[r]{-2}}
\put(-45,8){\makebox(10,5)[r]{-1}}
\put(-45,18){\makebox(10,5)[r]{0}}
\put(-45,28){\makebox(10,5)[r]{1}}
\put(-45,38){\makebox(10,5)[r]{2}}
\put(-45,48){\makebox(10,5)[r]{3}}
\put(-45,58){\makebox(10,5)[r]{4}}
\put(-45,68){\makebox(10,5)[r]{5}}
\put(-45,78){\makebox(10,5)[r]{6}}
\put(-45,88){\makebox(10,5)[r]{7}}
\put(-45,98){\makebox(10,5)[r]{8}}
\put(-45,108){\makebox(10,5)[r]{9}}
\put(-30,112){\line(0,1){1}}
\put(-30,114){\line(0,1){1}}
\put(-30,116){\line(0,1){1}}
\put(-30,-4){\line(0,1){1}}
\put(-30,-6){\line(0,1){1}}
\put(-30,-8){\line(0,1){1}}
\put(-30,0){\line(0,1){110}}
\multiput(-31,0)(0,10){12}{\line(1,0){2}}
\put(-1,40){\line(1,0){2}}
\put(0,40){\vector(0,1){70}}
\put(10,10){\vector(0,1){80}}
\put(9,10){\line(1,0){2}}
\put(20,40){\vector(0,1){50}}
\put(19,40){\line(1,0){2}}
\put(30,50){\vector(0,1){40}}
\put(29,50){\line(1,0){2}}
\put(40,50){\vector(0,1){40}}
\put(39,50){\line(1,0){2}}
\put(50,60){\vector(0,1){30}}
\put(49,60){\line(1,0){2}}
\put(60,70){\vector(0,1){20}}
\put(59,70){\line(1,0){2}}
\put(70,40){\vector(0,1){40}}
\put(69,40){\line(1,0){2}}
\put(80,50){\vector(0,1){30}}
\put(79,50){\line(1,0){2}}
\put(90,50){\vector(0,1){30}}
\put(89,50){\line(1,0){2}}
\put(100,70){\vector(0,1){10}}
\put(99,70){\line(1,0){2}}
\put(110,50){\vector(0,1){4}}
\put(109,50){\line(1,0){2}}
\put(120,0){\vector(0,1){40}}
\put(119,0){\line(1,0){2}}
\put(130,10){\vector(0,1){20}}
\put(129,10){\line(1,0){2}}
\end{picture}
\end{description}
If we say $\Del$ is in one of the orders above, we mean
its segments are weakly decreasing from left to right.

     Given a multisegment $\Del$, let
   $\xDel$ denote the multisegment in \lex order, and $\sDel$
   denote it in \std order.
We use \lex order when computing
$\eti$ and \std order when computing $\ethati$.

   Given a multisegment and a choice of ordering of its segments
   $\Del = \{\del{i_1}{j_1}, \ldots, $  
$ \del{i_m}{j_m} \}$, recall  that we let 
   $$\Ind \Del$$
   denote $\Ind_{n_1, n_2, \ldots, n_m}^n \del{i_1}{j_1} \bx \cdots
   \bx \del{i_m}{j_m}$,
   where $n_k = j_k - i_k +1$ and $n = \sum_k n_k$.  
Also $n(\Del) = n$ and $m(\Del) = m$.

   Observe that different orderings of a multisegment $\Del$
   can give non-isomorphic $\Ind \Del$.  For instance, let
   $\Del = \{ \del 00, \del 11 \}$.  We  know that 
$\Ind \sDel =$
   $\Ind_{1,1}^2 1 \bx q \not\isom \Ind_{1,1}^2 q \bx 1$
$= \Ind \xDel$,
   although they do (always) have the same composition factors.
   
The corollary below says
   The modules $\Ind \Del$ are isomorphic to each other for 
   the two orderings above.

   \begin{cor} 
   \label{cor-order}
   Let $\Del$ be a multisegment.  
Then
$\Ind \xDel \isom \Ind \sDel$.
   \end{cor}
   
   \begin{proof}
   From $\xDel$ in \lex order,
consider a segment $\del{i_0}{j_0}$ with largest
   \std 
   order, i.e.~ with largest $j$, and among those with smallest $i$.
   In terms of a picture, that means one (not necessarily unique)
   that reaches highest, and is longest among those.
(In the example above, that would be the eleventh segment $\del 29$.)
   Because $j_0$ is largest, all segments $\del{i}{j}$ 
   to the left 
   of $\del{i_0}{j_0}$  have $j \le j_0$; and because $\xDel$ is in \lex order, 
   we also have $i \ge i_0$.
   In other words $i_0 \le i \le j \le j_0$, so by part (ii) 
of the linking
  lemma $\Ind \del ij \bx \del{i_0}{j_0} \isom \Ind \del{i_0}{j_0} \bx \del ij$.
   By  transitivity   
   of induction, and repeating this argument,
the module $\Ind \xDel$ does not change if we slide
   $\del{i_0}{j_0}$ all the way over to the left, where it belongs with respect
   to \std order.
(One can see $\del 29$ is indeed leftmost in the example of
\std order above.)
   We can repeat this process with the remaining segments, eventually yielding
   $\Ind \xDel \isom \Ind \sDel$.
   \end{proof}
   
   
Theorem \ref{thm-class} stated in section \ref{sec-main}
shows that multisegments
parameterize the irreducible modules in $\RqHaffn$, $n \ge 0$.
In other words, 
given a multisegment $\Del$, the module $\Md = \cosoc \Ind \Del$
is irreducible, and as $\Del$ ranges over all multisegments,
$\Md$ ranges over all irreducibles in $\Rq$.

We will include the proof of theorem \ref{thm-class} here, starting
with the following lemma.
   
   \begin{lemma}
\label{lemma-j} 
Suppose $\Del = \{ \del{i_1}j,\del{i_2}j, \ldots, \del{i_m}j \}$.
   Then $\Ind \Del$ is irreducible.
   \end{lemma}
   \begin{proof}Let $\Delbar = \{\del{i_1}{j-1}, \ldots, \del{i_m}{j-1} \}$
   (entirely ommiting segments for which $i_k = j$).
   By induction, we may assume $\Ind \Delbar$ is irreducible.
Observe that the base case of this lemma is given by Kato's Theorem,
   which asserts $\Ind q^j \bx \cdots \bx q^j$ is irreducible.
   
   From the linking lemma,
 $\Ind \del{i_1}{j-1} \bx \del{i_2}j \bx \del jj 
   \isom  \Ind \del{i_1}{j-1} \bx \del jj \bx \del{i_2}j $
which surjects
   to $\Ind \del {i_1}j \bx \del{i_2} j$.
   Iterating this argument, we get a surjection
   $$\Ind \Delbar \bx q^j \bx q^j \bx \cdots \bx q^j \surj \Ind \Del \surj \cosoc \Ind \Del.$$
   Because $\Ind \Delbar$ is irreducible by the inductive
hypothesis, proposition \ref{lemma-3.5} 
shows that 
   $\Ind \Delbar \bx q^j  \bx \cdots \bx q^j $ has irreducible cosocle $M$. 
   Therefore we must also have that
   $\cosoc \Ind \Del = M$. 
   Furthermore,  $M$ occurs with multiplicity one and
   is the unique composition factor of
   $\Ind \Delbar \bx q^j  \bx \cdots \bx q^j$, and thus of
   $\Ind \Del$, which has $\epsj{j}(M) = m$ (also by 
proposition \ref{lemma-3.5}). 
   On the other hand, the linking lemma shows $\Ind \Del$ does
   not depend on the order of its segments, and so by 
   proposition \ref{prop-soccosoc}, $M = \soc \Ind \Del$ as well.
   Then by  the remark following 
proposition \ref{prop-irred}, $\Ind \Del = M$ is irreducible.
   \end{proof}

Observe that for $\Del$ as in the lemma above, $\epsj{j}(\Ind \Del) = m$,
and then by 
Corollary \ref{cor-e-et} and Lemma \ref{lemma-mackey},
$\et{j}^m (\Ind \Del) = \Ind \Delbar$.

Now we complete the proof of theorem \ref{thm-class}.
   \begin{proof}
   In light of corollary \ref{cor-order},
   we may prove the theorem for $\Md = \cosoc \Ind \sDel$,
   which we will do so  inducting on $n$ and making use of the functor $\ei$.
   
We will show $\Md$ is irreducible.
   First, we can write $\Del$ as a disjoint union of multisegments
   $\Del = \bigcup_j \Delsup{j}$ where $\Delsupj = \{ \del ik \in \Del
 \mid k = j \}$.
The lemma above showed $\Ind \Delsupj$ is irreducible.
   
   For $\Delsupj = \{ \del{i_1}j, \ldots, \del{i_m}j \}$ with
   $i_1 \le i_2 \le \cdots \le i_m$, define $\beta = \beta(\Delsupj)$ to be the 
   partition conjugate to $(j-i_1+1, j-i_2+1, \ldots, j-i_m+1)$.

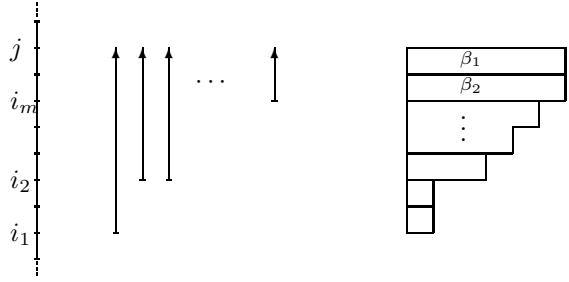
\begin{picture}(150,130)(-50,0)
\put(-30,112){\line(0,1){1}}  
\put(-30,114){\line(0,1){1}}
\put(-30,116){\line(0,1){1}}
\put(-30,14){\line(0,1){1}}
\put(-30,16){\line(0,1){1}}
\put(-30,18){\line(0,1){1}}
\put(-30,20){\line(0,1){90}}  
\put(-31,20){\line(1,0){2}}   
\put(-31,30){\line(1,0){2}}
\put(-31,40){\line(1,0){2}}
\put(-31,50){\line(1,0){2}}
\put(-31,60){\line(1,0){2}}
\put(-31,70){\line(1,0){2}}
\put(-31,80){\line(1,0){2}}
\put(-31,90){\line(1,0){2}}
\put(-31,100){\line(1,0){2}}
\put(-31,110){\line(1,0){2}}
\put(0,30){\vector(0,1){70}}
\put(-1,30){\line(1,0){2}}
\put(-40,27){$i_1$}
\put(-40,47){$i_2$}
\put(10,50){\vector(0,1){50}}
\put(9,50){\line(1,0){2}}
\put(-40,77){$i_m$}
\put(20,50){\vector(0,1){50}}
\put(19,50){\line(1,0){2}}
\put(30, 85){$\cdots$}
\put(60,80){\vector(0,1){20}}
\put(59,80){\line(1,0){2}}
\put(-40,97){$j$}
\put(110,100){\line(1,0){60}}
\put(110,30){\line(0,1){70}}
\put(110,90){\line(1,0){60}}
\put(170,90){\line(0,1){10}}
\put(110,80){\line(1,0){60}}
\put(170,80){\line(0,1){10}}
\put(130,94){{\scriptsize $\beta_1$}}
\put(130,84){{\scriptsize $\beta_2$}}
\put(130,65){$\vdots$}

\put(160,70){\line(0,1){10}}
\put(150,70){\line(1,0){10}}
\put(150,60){\line(0,1){10}}
\put(110,60){\line(1,0){40}}
\put(140,50){\line(0,1){10}}
\put(120,30){\line(0,1){20}}
\put(110,30){\line(1,0){10}}
\put(110,40){\line(1,0){10}}
\put(110,50){\line(1,0){30}}
\end{picture}

Let $\Qdel(\Delsupj)$ be the term
$ q^{i_1} q^{i_1+1} \cdots q^{i_2} q^{i_2}
   q^{i_2 +1} q^{i_2 +1} \cdots q^j q^j \cdots q^j$
and observe $\Qdel(\Delsupj)$ occurs
with multiplicity $\beta ! = \beta_1 ! \beta_2 !
\cdots \beta_{j-i_1+1}!$ 
   in $\ch \Ind \Delsupj$.
Also notice that the proof of the previous lemma, the observation made
following its proof, and induction show
\begin{gather}
\label{eqn-Q}
\1 = \Mg{\emptyset} = \et{i_1}^{\beta_{j-i_1+1}} \cdots
\et{j-1}^{\beta_2} \et{j}^{\beta_1} (\Ind \Delsupj). 
\end{gather}
   
   Now return to $\sDel = \Delsup{j_1} \cup \cdots \cup \Delsup{j_t}$
   with $j_1 > j_2 >\cdots > j_t$.  
   Because $j$ is strictly bigger than any exponent occurring in
the segments of  $\Delsup{j_a}$
   if $j > j_a$, from the shuffle lemma we see that the concatenation
   $\Qdel(\Delsup{j_1})\cdots \Qdel(\Delsup{j_t})$ occurs as a term in 
   $\ch \Ind \sDel$ with multiplicity $\beta(\Delsup{j_1})! \cdots \beta(\Delsup{j_t})!$.
This is its multiplicity in $\sDel$ as well.
   The irreducibility of $\Ind \Delsup{j_a}$ and Frobenius reciprocity show this
   term occurs in any quotient of $\Ind \sDel$ with that same multiplicity. 
   Therefore by proposition \ref{prop-irred}
 $\Md = \cosoc \Ind \sDel = \cosoc \Ind \xDel$
must be irreducible (and furthermore occurs with
multiplicity one in $\Ind \sDel$).

For parts 2 and 3 of the theorem, we give only a brief sketch
of the argument from \cite{BZ,Z}.  We will see later that it also
follows from the computation for $\eti \Md$.  

Next, we will show irreducibles corresponding to distinct multisegments
are distinct.
Given an irreducible module $\Md$, let $j$ be the {\it  smallest\/}
 integer for which $\epsj{j}(\Md) >0$. 
Write $\ve = \epsj{j}(\Md)$. Let $N = \et{j}^{\ve} \Md$.
From lemma \ref{lemma-j} and the fact that we may assume $\Del$ is
in \std order, we must have 
$\ve = m(\Delsupj) = m$.
Using lemma \ref{lemma-mackey} and corollary \ref{cor-e-et}, we
must have that $N = \Mg{\Gamma}$ where
$$\Gamma =  (\Delsupj)^- \cup
\bigcup_{k \neq j} \Delsup{k}.$$
Notice all the $k$ we union over are bigger than $j$ by our choice of $j$.
Suppose $\Md \isom \Mg{\Del'}$.  
Then also $N =   \et{j}^{\ve} \Mg{\Del'} = \Mg{\Gamma'}$
for $\Gamma'$ obtained from $\Del'$ in the same way as above.
By induction $\Gamma = \Gamma'$.
Lemma \ref{lemma-j} gives us
a well-defined and reversible
 way to go from $\Del$ to $\Gamma$.
By the properties of \std order,
if $\Gamma$ contains $k$ segments of the form $\del i{j-1}$,
each one must have come from a $\del ij$ in $\Del$, and further
$\Del$ must contain $\ve - k$ segments $\del jj$.
This holds as well for $\Del'$, so
that $\Md \isom \Mg{\Del'} \implies \Del = \Del'$.

Finally, for part 3, we will show every irreducible module in $\Rq$ is some
$\Md$.
Let $M$ be any irreducible module in $\Rq$, and let 
$N = \et{j}^{\epsj{j}(M)} M$, where again $j$ is smallest possible.
By induction there is some multisegment $\Gamma$
such that  $N = \Mg{\Gamma}$.

We construct a new multisegment $\Del$ from $\Gamma$ by replacing
each $\del i{j-1} \in \Gamma$ by $\del ij$, and also adding in
$\epsj{j}(M) - \epsj{j-1}(N)$ many segments  $\del jj$.  
This number is nonnegative because  $\epsj{j-1}(M) = 0$ by choice of $j$
and proposition \ref{prop-10.4} shows that
   $\epsj{j-1}(M) + k \ge \epsj{j-1}(\et{j}^k(M))$.
Then $\Ind \sGamma \bx q^j \bx \cdots \bx q^j$ has cosocle $M$,
but
by an argument similar to that in the proof of lemma \ref{lemma-j},
it also surjects to $\Ind \sDel$,
which has cosocle $\Md$.  Thus $M = \Md$.
   \end{proof}
   
   \begin{cor}
\label{cor-multone}
$\Md$ occurs with multiplicity one as a composition factor of $\Ind \Del$.
   \end{cor}
   
\begin{rem}
\label{rem-j}
Because $\Md = \cosoc \Ind \sDel$, 
if $\del ij \in \Del$ is smallest in  
 \std order, then $j$ 
is the smallest integer such that $\epsj{j}(\Md) \neq 0$ and furthermore
\begin{gather*}
\epsj{j}(\Md) = \mid  \{ \del ij \in \Del \mid i \in \Z \} \mid
= m(\Delsupj).
\end{gather*}
For other $j$ the cardinality of this set is merely an upper bound
on $\epsj{j}(\Md)$.
\end{rem}

In the proof of parts 2 and 3 of the above theorem, we saw that
lemma \ref{lemma-j} allowed us to compute the multisegment
$\Gamma$ for which  $\et{j}^{\epsj{j}(\Md)} \Md = \Mg{\Gamma}$,
where  $j$ was the {\it smallest\/}  integer for which $\epsj{j}(\Md) >0$.
However, at this stage
we have not proved that Rule \ref{defn-Ei} describes $\eti \Md$
in terms of multisegments for $i \neq j$, nor even
$\et{j}^k \Md$ when $k \neq \epsj{j}(\Md)$.   
This will be proved below. 
However having done so, it will be clear theorem \ref{thm-Ei} 
also gives a proof of parts 2 and 3
of theorem \ref{thm-class}.

\begin{rem}
\label{rem-crystal}
One could now appeal to
Grojnowski's Theorem 14.3: that the graph 
whose nodes are indexed by multisegments $\Del$ and whose edges
are given by $\Gamma \xrightarrow{i} \Del$ if
$$\et{i} \Md = \Mg{\Gamma}$$
is the crystal graph $\B(\infty)$.

Given a multisegment $\Del$, we can iterate the rule that we always
take $\et{j}^{\epsj{j}}$ for the smallest $j$ such that
$\epsj{j} \neq 0$. 
Combinatorially, we are just iterating the replacement of
$\Delsup{j}$ with $(\Delsup{j})^-$,  (if we wrote
$\Del = \bigcup_k \Delsup{k}$).
This   constructs a distinguished
path on the crystal graph from $\Del$ back to the empty multisegment
$\emptyset$. 
Now that we have determined where on the crystal graph the node
$\Del$ sits, the determination of a single edge 
$$ ? \xrightarrow{i} \Del$$
is purely combinatorial.
And conversely, we can determine the distinguished path leading
from $?$ back to $\emptyset$, and thus re-express $\Ei \Del$
as a multisegment.
This is exactly given by rule \ref{defn-Ei}.
\end{rem}

  Theorem  \ref{thm-Ei} and rule \ref{defn-Ei} give
the combinatoric that Grojnowski's theorem
dictates, but the proof here is
module-theoretic.
 Admittedly, a proof appealing to  the known crystal structure
described in Theorem 14.3 of \cite{G} is much slicker.

\section{Example of computing $\Ei$}
\label{sec-ex}
Here we give an example of computing $\Ei \Del$.

Let $\Del =
\{ \del 56, \del 57, \del 47, \del 33, \del 36, \del 36, \del 37,
\del 37, \del 26, \del 27, \del 29,$ $\del {-1}7,$
$ \del {-1}1, \del{-2}2 \}$
as in the example for \lex order in section \ref{sec-class}.

\begin{picture}(180,140)(-110,-20)
\put(-60,112){\line(0,1){1}}
\put(-60,114){\line(0,1){1}}
\put(-60,116){\line(0,1){1}}
\put(-60,-4){\line(0,1){1}}
\put(-60,-6){\line(0,1){1}}
\put(-60,-8){\line(0,1){1}}
\put(-60,0){\line(0,1){110}}
\put(-61,0){\line(1,0){2}}
\put(-61,10){\line(1,0){2}}
\put(-61,20){\line(1,0){2}}
\put(-61,30){\line(1,0){2}}
\put(-61,40){\line(1,0){2}}
\put(-61,50){\line(1,0){2}}
\put(-61,60){\line(1,0){2}}
\put(-61,70){\line(1,0){2}}
\put(-61,80){\line(1,0){2}}
\put(-61,90){\line(1,0){2}}
\put(-61,100){\line(1,0){2}}
\put(-61,110){\line(1,0){2}}
\put(-75,-2){\makebox(10,5)[r]{-2}}
\put(-75,8){\makebox(10,5)[r]{-1}}
\put(-75,18){\makebox(10,5)[r]{0}}
\put(-75,28){\makebox(10,5)[r]{1}}
\put(-75,38){\makebox(10,5)[r]{2}}
\put(-75,48){\makebox(10,5)[r]{3}}
\put(-75,58){\makebox(10,5)[r]{4}}
\put(-75,68){\makebox(10,5)[r]{5}}
\put(-75,78){\makebox(10,5)[r]{6}}
\put(-75,88){\makebox(10,5)[r]{7}}
\put(-75,98){\makebox(10,5)[r]{8}}
\put(-75,108){\makebox(10,5)[r]{9}}
\put(-115,48){\makebox(10,15)[l]{$\Del$}}
\put(-32, -5){$+$}
\put(-22, -5){$-$}
\put(-12, -5){$-$}
\put(8, -5){$+$}
\put(18, -5){$+$}
\put(28, -5){$-$}  
\put(38, -5){$-$}
\put(48, -5){$+$}
\put(58, -5){$-$}  
\put(78, -5){$-$}  
\put(-32, -15){$+$}
\put(28, -15){$-$}  
\put(58, -15){$-$}  
\put(78, -15){$-$}  
\put(-31,70){\line(1,0){2}}
\put(-30,70){\vector(0,1){10}}
\put(-21,70){\line(1,0){2}}
\put(-20,70){\vector(0,1){20}}
\put(-11,60){\line(1,0){2}}
\put(-10,60){\vector(0,1){30}}
\put(-1,50){\line(1,0){2}}
\put(0,50){\vector(0,1){4}}
\put(10,50){\vector(0,1){30}}
\put(9,50){\line(1,0){2}}
\put(20,50){\vector(0,1){30}}
\put(19,50){\line(1,0){2}}
\put(30,50){\vector(0,1){40}}  
\put(29,50){\line(1,0){2}}
\put(40,50){\vector(0,1){40}}
\put(39,50){\line(1,0){2}}
\put(50,40){\vector(0,1){40}}
\put(49,40){\line(1,0){2}}
\put(60,40){\vector(0,1){50}}  
\put(59,40){\line(1,0){2}}
\put(70,40){\vector(0,1){70}}
\put(69,40){\line(1,0){2}}
\put(80,10){\vector(0,1){80}}  
\put(79,10){\line(1,0){2}}
\put(90,10){\vector(0,1){20}}
\put(89,10){\line(1,0){2}}
\put(100,0){\vector(0,1){40}}
\put(99,0){\line(1,0){2}}
\end{picture}

We will compute $\E{7}^k \Del$.
First, we compute the word
$+-- \; ++--+- \; - \; \;$ and cancel
as $+(-(- \; +)+)-(-+)- \; - \; \;$
leaving $+ \; \; \;  \;- \; \;- \; - \; \;$ uncanceled.

Thus 
$\epsj{7}(\Md) = 3$ and 
 $\E{7}^4 \Del = \0$.
%

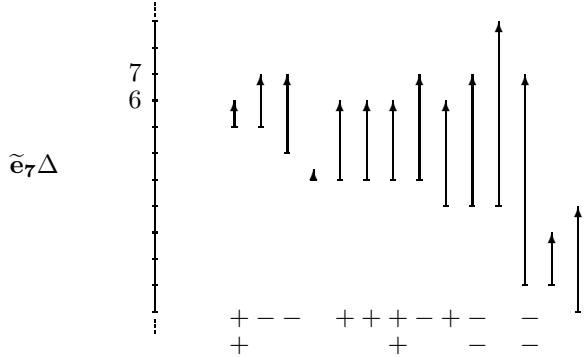
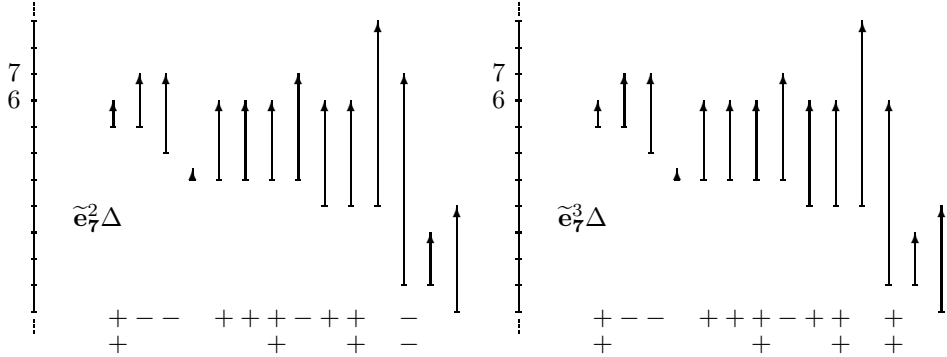
\begin{picture}(180,140)(-110,-20)
\put(-60,112){\line(0,1){1}}
\put(-60,114){\line(0,1){1}}
\put(-60,116){\line(0,1){1}}
\put(-60,-4){\line(0,1){1}}
\put(-60,-6){\line(0,1){1}}
\put(-60,-8){\line(0,1){1}}
\put(-60,0){\line(0,1){110}}
\put(-61,0){\line(1,0){2}}
\put(-61,10){\line(1,0){2}}
\put(-61,20){\line(1,0){2}}
\put(-61,30){\line(1,0){2}}
\put(-61,40){\line(1,0){2}}
\put(-61,50){\line(1,0){2}}
\put(-61,60){\line(1,0){2}}
\put(-61,70){\line(1,0){2}}
\put(-61,80){\line(1,0){2}}
\put(-61,90){\line(1,0){2}}
\put(-61,100){\line(1,0){2}}
\put(-61,110){\line(1,0){2}}
\put(-75,78){\makebox(10,5)[r]{6}}
\put(-75,88){\makebox(10,5)[r]{7}}
\put(-115,48){\makebox(10,15)[l]{$\E{7} \Del$}}
\put(-32, -5){$+$}
\put(-22, -5){$-$}
\put(-12, -5){$-$}
\put(8, -5){$+$}
\put(18, -5){$+$}
\put(28, -5){$+$}  
\put(38, -5){$-$}
\put(48, -5){$+$}
\put(58, -5){$-$}
\put(78, -5){$-$}
\put(-32, -15){$+$}
\put(28, -15){$+$}  
\put(58, -15){$-$}
\put(78, -15){$-$}
\put(-31,70){\line(1,0){2}}
\put(-30,70){\vector(0,1){10}}
\put(-21,70){\line(1,0){2}}
\put(-20,70){\vector(0,1){20}}
\put(-11,60){\line(1,0){2}}
\put(-10,60){\vector(0,1){30}}
\put(-1,50){\line(1,0){2}}
\put(0,50){\vector(0,1){4}}
\put(10,50){\vector(0,1){30}}
\put(9,50){\line(1,0){2}}
\put(20,50){\vector(0,1){30}}
\put(19,50){\line(1,0){2}}
\put(30,50){\vector(0,1){30}}  
\put(29,50){\line(1,0){2}}
\put(40,50){\vector(0,1){40}}
\put(39,50){\line(1,0){2}}
\put(50,40){\vector(0,1){40}}
\put(49,40){\line(1,0){2}}
\put(60,40){\vector(0,1){50}}
\put(59,40){\line(1,0){2}}
\put(70,40){\vector(0,1){70}}
\put(69,40){\line(1,0){2}}
\put(80,10){\vector(0,1){80}}
\put(79,10){\line(1,0){2}}
\put(90,10){\vector(0,1){20}}
\put(89,10){\line(1,0){2}}
\put(100,0){\vector(0,1){40}}
\put(99,0){\line(1,0){2}}
\end{picture}


\begin{picture}(160,150)(-60,-20)
\put(-60,112){\line(0,1){1}}
\put(-60,114){\line(0,1){1}}
\put(-60,116){\line(0,1){1}}
\put(-60,-4){\line(0,1){1}}
\put(-60,-6){\line(0,1){1}}
\put(-60,-8){\line(0,1){1}}
\put(-60,0){\line(0,1){110}}
\put(-61,0){\line(1,0){2}}
\put(-61,10){\line(1,0){2}}
\put(-61,20){\line(1,0){2}}
\put(-61,30){\line(1,0){2}}
\put(-61,40){\line(1,0){2}}
\put(-61,50){\line(1,0){2}}
\put(-61,60){\line(1,0){2}}
\put(-61,70){\line(1,0){2}}
\put(-61,80){\line(1,0){2}}
\put(-61,90){\line(1,0){2}}
\put(-61,100){\line(1,0){2}}
\put(-61,110){\line(1,0){2}}
\put(-75,78){\makebox(10,5)[r]{6}}
\put(-75,88){\makebox(10,5)[r]{7}}
\put(-45,28){\makebox(10,15)[l]{$\E{7}^2 \Del$}}
\put(-32, -5){$+$} 
\put(-22, -5){$-$} 
\put(-12, -5){$-$}
\put(8, -5){$+$} 
\put(18, -5){$+$} 
\put(28, -5){$+$}  
\put(38, -5){$-$} 
\put(48, -5){$+$} 
\put(58, -5){$+$}  
\put(78, -5){$-$}
\put(-32, -15){$+$} 
\put(28, -15){$+$}  
\put(58, -15){$+$}  
\put(78, -15){$-$}
\put(-31,70){\line(1,0){2}}
\put(-30,70){\vector(0,1){10}}
\put(-21,70){\line(1,0){2}}
\put(-20,70){\vector(0,1){20}}
\put(-11,60){\line(1,0){2}}
\put(-10,60){\vector(0,1){30}}
\put(-1,50){\line(1,0){2}}
\put(0,50){\vector(0,1){4}}
\put(10,50){\vector(0,1){30}}
\put(9,50){\line(1,0){2}}
\put(20,50){\vector(0,1){30}}
\put(19,50){\line(1,0){2}}
\put(30,50){\vector(0,1){30}}  
\put(29,50){\line(1,0){2}}
\put(40,50){\vector(0,1){40}}
\put(39,50){\line(1,0){2}}
\put(50,40){\vector(0,1){40}}
\put(49,40){\line(1,0){2}}
\put(60,40){\vector(0,1){40}}  
\put(59,40){\line(1,0){2}}
\put(70,40){\vector(0,1){70}}
\put(69,40){\line(1,0){2}}
\put(80,10){\vector(0,1){80}}
\put(79,10){\line(1,0){2}}
\put(90,10){\vector(0,1){20}}
\put(89,10){\line(1,0){2}}
\put(100,0){\vector(0,1){40}}
\put(99,0){\line(1,0){2}}
\end{picture}
%
\begin{picture}(160,150)(-80,-20)
\put(-60,112){\line(0,1){1}}
\put(-60,114){\line(0,1){1}}
\put(-60,116){\line(0,1){1}}
\put(-60,-4){\line(0,1){1}}
\put(-60,-6){\line(0,1){1}}
\put(-60,-8){\line(0,1){1}}
\put(-60,0){\line(0,1){110}}
\put(-61,0){\line(1,0){2}}
\put(-61,10){\line(1,0){2}}
\put(-61,20){\line(1,0){2}}
\put(-61,30){\line(1,0){2}}
\put(-61,40){\line(1,0){2}}
\put(-61,50){\line(1,0){2}}
\put(-61,60){\line(1,0){2}}
\put(-61,70){\line(1,0){2}}
\put(-61,80){\line(1,0){2}}
\put(-61,90){\line(1,0){2}}
\put(-61,100){\line(1,0){2}}
\put(-61,110){\line(1,0){2}}
\put(-75,78){\makebox(10,5)[r]{6}}
\put(-75,88){\makebox(10,5)[r]{7}}
\put(-45,28){\makebox(10,15)[l]{$\E{7}^3 \Del$}}
\put(-32, -5){$+$} 
\put(-22, -5){$-$} 
\put(-12, -5){$-$}
\put(8, -5){$+$} 
\put(18, -5){$+$} 
\put(28, -5){$+$}  
\put(38, -5){$-$} 
\put(48, -5){$+$} 
\put(58, -5){$+$}  
\put(78, -5){$+$}  
\put(-32, -15){$+$} 
\put(28, -15){$+$}  
\put(58, -15){$+$}  
\put(78, -15){$+$}  
\put(-31,70){\line(1,0){2}}
\put(-30,70){\vector(0,1){10}}
\put(-21,70){\line(1,0){2}}
\put(-20,70){\vector(0,1){20}}
\put(-11,60){\line(1,0){2}}
\put(-10,60){\vector(0,1){30}}
\put(-1,50){\line(1,0){2}}
\put(0,50){\vector(0,1){4}}
\put(10,50){\vector(0,1){30}}
\put(9,50){\line(1,0){2}}
\put(20,50){\vector(0,1){30}}
\put(19,50){\line(1,0){2}}
\put(30,50){\vector(0,1){30}}  
\put(29,50){\line(1,0){2}}
\put(40,50){\vector(0,1){40}}
\put(39,50){\line(1,0){2}}
\put(50,40){\vector(0,1){40}}
\put(49,40){\line(1,0){2}}
\put(60,40){\vector(0,1){40}}  
\put(59,40){\line(1,0){2}}
\put(70,40){\vector(0,1){70}}
\put(69,40){\line(1,0){2}}
\put(80,10){\vector(0,1){70}}  
\put(79,10){\line(1,0){2}}
\put(90,10){\vector(0,1){20}}
\put(89,10){\line(1,0){2}}
\put(100,0){\vector(0,1){40}}
\put(99,0){\line(1,0){2}}
\end{picture}

As an exercise, the reader can consider $\sDel$ and compute
$\Ehat{-1}(\Del) = 2$,
$\Ehat{2} = 2$,
$\Ehat{3} = 4$,
$\Ehat{4} = 1$,
$\Ehat{5} = 2$.


\section{Proof of theorem \ref{thm-Ei}}
\label{sec-proofseg}
   \begin{rem}
\label{rem-cosoc}
In this section and the next, we will repeatedly use the fact that
if $\cosoc \Ind A \bx B$ is irreducible, then it coincides
with $\cosoc (\Ind (\cosoc A) \bx B)$ and $\cosoc (\Ind A \bx (\cosoc B))$.
   \end{rem}
   
   \begin{lemma}[sliding lemma]
\label{lemma-triple}
\begin{enumerate}
\item
Let $\Del = \{ \del bj, \del a{j-1} \}$.
Then if $a < b \le c$ 
$$\Ind \Md \bx \del cj \isom \Ind \del cj \bx \Md$$
is irreducible.
\item
Let $\Gamma = \{ \del a{z_0}, \del {a-1}{z_1}, \ldots, \del {a-k}{z_k} \}$,
where $z_0 > z_1 > \cdots > z_k \ge z$. Then
$$\Ind \Mg{\Gamma} \bx \del{a-k}z \isom \Ind  \del{a-k}z \bx\Mg{\Gamma}$$
is irreducible.
\end{enumerate}
   \end{lemma}
   
   \begin{proof}
Write $M = \Md$.  Let $N = \cosoc \Ind \del cj \bx M$.
Then $\del cj, \del bj, \del a{j-1}$ are in \lex order, so
that by the remark \ref{rem-cosoc}
 above, $N = \Mg{\{\del cj, \del bj, \del a{j-1}\}}$.

First, we claim $\epsj{j} M = 0$ and  $\epsj{j} N = 1$.
We have an exact sequence
$$ 0 \to \Ind \del b{j-1} \bx \del a{j} \to \Ind \del b{j} \bx \del a{j-1} 
\to M \to 0$$
by the linking lemma.
Applying the exact functor $\e{j}$ and using lemma \ref{lemma-mackey}
yields
$$ 0 \to \Ind \del b{j-1} \bx \del a{j-1} \to \Ind \del b{j-1} \bx \del a{j-1} 
\to \e{j}M \to 0.$$
Because the first two terms are   irreducible and isomorphic by
the linking lemma, the
last term must be zero and thus $\epsj{j} M = 0$.

By the shuffle lemma  $\epsj{j} N \le 1$. On the other hand,
since $\Ind \del cj \bx \del bj = \Ind \del bj \bx \del cj$, we
know by proposition \ref{prop-ch} that $\ch N$ has a term 
$q^b q^{b+1} \cdots q^c q^c q^{c+1} q^{c+1} \cdots q^j q^j q^a \cdots q^{j-1}$
and therefore a term $q^b \cdots q^j q^j q^{j-1}$, so finally
$q^b \cdots q^j q^{j-1} q^j$  
by remark \ref{rem-serre}. Thus    $\epsj{j} N \ge 1$.

Because $\epsj{j} N = 1$, $\e{j} N = \et{j} N$ is irreducible (and nonzero)
and we have by lemma \ref{lemma-mackey} 
$$\Ind \del c{j-1} \bx M \surj \e{j} N = \et{j} N,$$
which also shows that $\et{j} N =
\Mg{\{\del c{j-1}, \del bj \del a{j-1}\}}$ since these segments
are in \lex order.

Let $L = \cosoc \Ind M \bx \del cj$.
Then as above we get
$$\Ind M \bx \del c{j-1} \surj \e{j} L = \et{j} L.$$
Now $\del bj, \del a{j-1}, \del c{j-1}$ are in \std order,
so by theorem \ref{thm-class} and corollary \ref{cor-order},
$\et{j} L \isom \Mg{ \{ \del bj, \del a{j-1}, \del c{j-1} \} }
\isom \Mg{ \{ \del c{j-1},  \del bj, \del a{j-1} \} }
\isom  \et{j} N$. (In particular, this shows
$L$ was irreducible, since any summands of it must
have  $\epsj{j}= 1$ as well.)  Therefore, from proposition
\ref{cor-3.6} 
we know $L \isom N$.

On the other hand, $\cosoc \Ind \del cj \bx M = N
\isom L = \soc \Ind \del cj \bx M$ by proposition \ref{prop-soccosoc}.
However by corollary \ref{cor-multone}, $N$ occurs with
multiplicity one in $\Ind \del c{j} \bx \del bj  \bx \del a{j-1}$
and therefore in $\Ind \del cj \bx M$.  From remark following
proposition
\ref{prop-irred},
$\Ind \del cj \bx M$ is irreducible, and thus also isomorphic
to $\Ind M \bx \del cj$.

This proves   part 1.  The proof of part 2 is quite similar, but we
use $\ehati$ instead of $\ei$.  Notationally, it is slightly
messier, but it is more useful  in this form  for use in  future theorems.

Now renotate, letting 
$M = \Mg{\Gamma}$ and let $N = \cosoc \Ind \Mg{\Gamma} \bx \del{a-k}{z}$.
Write $\Gamma' = \{ \del a{z_0},  \ldots, \del {a-k}{z_k}, \del{a-k}{z} \}$.
Then $\Gamma'$ is in \std order, so that $N = \Mg{\Gamma'}$ is irreducible.

The same argument as in part 1 shows $\epshat{a-i}(\Mg{\Gamma})= 0$ if $i > 0$,
and so $\epshat{a-k}(N) \le 1$.
An argument similar to that in part 1 also shows that
$\Ind \Mg{\{\del{a-k+1}{z_{k-1}}, \del{a-k}{z_k} \}} \bx \del{a-k}{z}
\isom \Ind \del{a-k}{z} \bx\Mg{\{\del{a-k+1}{z_{k-1}}, \del{a-k}{z_k} \}}$.
Together with remark \ref{rem-serre} and proposition \ref{prop-ch},
this shows $\epshat{a-k}(N) = 1$.

Let $L = \cosoc \Ind\del{a-k}{z} \bx \Mg{\Gamma}$. 
Observe that
$\Ind \Mg{\Gamma} \bx \del az = \Ind \del az \bx \Mg{\Gamma} $
by corollary \ref{cor-order} if $a \le z$. 
Let $b = \min \{a-1,z\}$, so that 
$\ehat{b} \ehat{b+1} \cdots \ehat{a-k} \del{a-k}{z} = \del a z$,
where we interpret $\del az = \triangle_\emptyset = \1$ if $a >z$. 
(If $a > z$, then $b = z$ and interpreting
$\del az = \1$, the observation still holds.)  By repeated application
of the exact functor $\ehati$ and the analogue to lemma \ref{lemma-mackey},
the cosocle of the above module is equal to each of the following: 
$\ethat{b} \cdots   \ethat{a-k+1}    \ethat{a-k} L
= \ehat{b} \cdots    \ehat{a-k} L$
$ = \cosoc \Ind \del az \bx \Mg{\Gamma} $
$= \ehat{b} \cdots    \ehat{a-k} N
= \ethat{b} \cdots    \ethat{a-k} N$.
That the right hand side is nonzero is similar to the 
computation that $\epshat{a-k}(N) \ge 1$. 
%
Proposition \ref{cor-3.6} implies that $L \isom N$.
A repeat of the argument that concludes part 1 shows that 
$\Ind \Mg{\Gamma} \bx \del{a-k}z \isom \Ind  \del{a-k}z \bx\Mg{\Gamma}$
is irreducible.
   \end{proof}
   
The previous lemma is the final step toward computing
 $\et{j} \Md$ 
where $\Del = \Delsupj \cup \Delsup{j-1}$, after which we can compute
$\et{j} \Md$ for
any $\Md$. 
Then we shall show it coincides with rule \ref{defn-Ei}.


\omitt{this can be explained by paths in perfect crystals...}

  \begin{lemma}
\label{lemma-ej}
For $\Del = \Delsupj \cup \Delsup{j-1},$
$$\et{j} \Md = \Mg{\E{j} \Del}$$
and $\et{j} \Md = 0$ if $\E{j} \Del =\0$.
In particular, $\epsj{j}(\Md) $ is the number of
$-$ signs left after all  cancelling is done as in rule \ref{defn-Ei}.
   \end{lemma}
  
   \begin{proof}
First we prove the second statement.
We have $\Md = \cosoc \Ind \xDel = \cosoc \Ind \del{a_1}{\pm}
\bx \cdots \bx \del{a_m}{\pm}$,
writing $-$ for $j$ and $+$ for $j-1$.
The $a_{i_k}$ weakly decrease, since we take $\Del$ in \lex order.
Recall from the proof of the sliding lemma \ref{lemma-triple}
that $\epsj{j}( \Mg{ \{\del{a_{i-1}}{j}, \del{a_i}{j-1} \}}) = 0$
if $a_{i-1} > a_i$,  in other words, if the pair is linked.
So, using remark \ref{rem-cosoc}, whenever you see an adjacent
$\del{a_{i-1}}{-} \bx  \del{a_i}{+}$, it does not contribute
to $\epsj{j}( \Md)$.
In other words, we can replace each $-+$ by $M_i = 
\Mg{ \{\del{a_{i-1}}{j}, \del{a_i}{j-1} \}}$, yielding
$\Md = \cosoc \Ind \del{a_1}{\pm} \bx \cdots \bx M_i \bx \cdots 
\del{a_m}{\pm}$.  Repeat.
Then by the sliding lemma, a  $\del{a_{i-k}}{-}$ to the left
of  $M_i$ can {\it slide\/} past it, leaving the cosocle of the
induced module unchanged.

Then we may repeat the first step, which in effect cancels
any adjacent $-+$ pairs, and then  by the sliding lemma
continues to cancel newly
created $-+$ pairs from the remaining symbols, and so on.
In the end we are left with something that looks like 
\begin{gather}
\label{eqn-form}
\Md  = \cosoc \Ind 
\underbrace{\del{a_1}{+} \cdots \bx M_i \cdots \bx \del{a_k}{+} }_{A}
\bx \underbrace{  \del{a_{k+1}}{-}
\bx \cdots  \bx \del{a_m}{-} }_{B} 
\end{gather}
where $A$ is an induced module consisting solely
of $\del{a_i}{+}$'s and $M_i$'s (which have the form
$\Mg{ \{ \del {a_{i-k}}{-}, \del {a_i}{+} \} }$), and $B$ is induced from
$\del{a_i}{-}$'s.
Thus $\epsj{j}(A) = 0$ and $\epsj{j}(B) = \epsj{j}(\Md)$ is the number
of uncanceled $-$ signs remaining from the procedure.
If there are no uncanceled $-$ signs left, then it is clear
$\epsj{j}(\Md) = 0$, so that $\et{j} \Md = 0 = \Mg{\0}= \Mg{\E{j} \Del}$.
This proves the second statement.
Notice that it was very important that $\xDel$ was in
\lex order (and not \std order) to start with.

However,  we have also shown that if we alter \lex order by sliding
segments from left to right past linked pairs, uncanceled $\pm$ symbols
will be unchanged.  


Now we can determine $\et{j} \Md$.  Let $\ve = \epsj{j}(\Md)$.
We know $\Md = \cosoc \Ind A \bx \del{a_{i_1}}j \bx \cdots 
\bx  \del{a_{i_\ve}}j$, and $\epsj{j}(A) = 0$.
Lemma \ref{lemma-mackey}
 tells us that $\e{j}(\Ind  A \bx \del{a_{i_1}}j \bx \cdots 
\bx  \del{a_{i_\ve}}j)$ is filtered by modules whose successive
quotients
have the form
\begin{gather}
\label{eqn-filter}
\Ind A \bx  \del{a_{i_1}}j \bx \cdots 
\bx \del{a_{i_s}}{j-1} \bx \cdots \bx \del{a_{i_\ve}}j,
\end{gather}
for  $s = 1, \ldots, \ve$.
Thus $\et{j} \Md$ is in the cosocle of one of these.
\omitt{Observe the $a_{i_k}$ are still in decreasing order.}
We will show below that each of these induced modules of the form
\eqref{eqn-filter} has irreducible
cosocle and so the above procedure lets us compute $\epsj{j}$ of those
cosocles.

We will show the module in \eqref{eqn-filter} is a quotient
of $\Ind \Del'$, where
 $\Del' = \Del  \cup \{ \del{a_{i_s}}{j-1}\} \setminus
	\{\del{a_{i_s}}{j} \}$ (in \lex order),
and so its cosocle is $\Mg{\Del'}$, which is irreducible.
Then we may apply the first part of the proof to compute
$\epsj{j}(\Mg{\Del'})$.

Starting with $\Ind \Del'$, which we know has irreducible 
cosocle, we will apply the same moves we performed on
$\Del$, i.e.~ ignore $-+$ pairs by forming $M_i$ to end
with a module as in \eqref{eqn-filter} as a quotient of 
$\Ind \Del'$.
We justify being able to perform these moves as follows.
First, by the linking lemma,
 $\Ind \Del'$ is the same whether $\Del'$ is in \lex order,
or in the order inherited from $\Del$ in \lex order with
$\del{a_{i_s}}{j}$  replaced by $\del{a_{i_s}}{j-1}$.
So let us start with $\Del'$ in the latter order.
Observe that when $c < b \le a$,
two applications of the linking lemma yield
$\Ind \del bj \bx \del{c}{j-1} \bx \del{a}{j-1} \isom
\Ind \del{a}{j-1} \bx \del bj \bx \del{c}{j-1}$,
the first of which is in \std order and so has irreducible
cosocle.
Therefore by remark \ref{rem-cosoc}, 
$$\cosoc \Ind \Mg{ \{ \del bj, \del{c}{j-1} \} } \bx \del{a}{j-1}
\isom \cosoc \Ind \del{a}{j-1} \bx \Mg{ \{ \del bj, \del{c}{j-1} \} }.
$$
If we compare this to the sliding lemma (lemma \ref{lemma-triple}), it says
any time we slid  $\del {a}j$ past  an $M_{i}$ in the above
process for $\Del$ (cancelling a $-+$), we would also have been allowed
to slide the segment past if it were replaced by $\del{a}{j-1}$.
So, starting from $\Del'$ we form the same $M_i$ and do the
same slidings to end with 
the module in equation \eqref{eqn-filter} as a quotient of $\Ind \Del'$.
Hence
$$\Mg{\Del'} = \cosoc \Ind A \bx \del{a_{i_1}}j \bx \cdots \bx
\del{a_{i_s}}{j-1} \bx\cdots\bx
\del{a_{i_\ve}}{j}.$$
In particular, the cosocle is irreducible.

  Now we may apply what we have already proved to compute
$\epsj{j}(\Mg{\Del'}) \le \ve - 2$ if $s \neq 1$.


However, we know from proposition \ref{prop-10.4} that 
$\epsj{j}(\et{j} \Md) = \ve -1$, so we must have $\et{j} \Md$ 
occurring only as
the cosocle of the induced module as in \eqref{eqn-filter} with  $s=1$.
In other words, it is the  {\it leftmost\/} remaining (uncanceled)
$-$ sign that signifies which segment to alter in computing $\E{j} \Del$. 
This process matches rule \ref{defn-Ei}, and so proves
the lemma.
(The only difference in \eqref{eqn-form} and rule \ref{defn-Ei}
 is that we slid $-+$ (i.e.~linked) pairs from the right to the left
of some  $-$ signs.
However it is the uncanceled symbols and the segments they correspond to
that are used in both computations, identically.)
   \end{proof}

We will now use lemma \ref{lemma-ej} to complete the proof of 
theorem \ref{thm-Ei}, that  $\et{j} \Md = \Mg{\E{j} \Del}$.
   \begin{proof}
From corollary \ref{cor-order}, we know
$\Md = \cosoc \Ind \xDel = \cosoc \Ind \sDel
= \cosoc \Ind \Delsup{k_1} \bx \cdots \bx \Delsupj \bx \Delsup{j-1} \bx
\cdots \bx \Delsup{k_s}$.
Let $j = k_t, j-1 = k_{t+1}$.
Let $A = \cosoc  \Ind \Delsup{k_1} \bx \cdots \bx \Delsup{k_{t-1}}$
and $B =\cosoc  \Ind \Delsup{k_{t+2}} \bx \cdots \bx \Delsup{k_s}$,
so that by remark \ref{rem-cosoc}
\begin{gather*}
\Md = \cosoc \Ind A \bx \Delsupj \bx \Delsup{j-1} \bx B.
\end{gather*}
Clearly $\epsj{j}(A) = 0$ and $\epsj{j}(B) = 0$ by the shuffle lemma.
Further, all terms in $\ch B$ have support consisting of 
$q^k$ where $k < j-1$.
Then remark \ref{rem-serre} implies
$\Ind q^j \bx  B \isom \Ind   B\bx  q^j$ (which furthermore is irreducible).
In other words,
$\epsj{j}(\Md) = \epsj{j}(\cosoc \Ind \Delsupj \bx \Delsup{j-1})$
and we have reduced theorem \ref{thm-Ei} to lemma \ref{lemma-ej}, where
the same argument holds.
To compute $\et{j}$ we now consider $\Delsupj \cup \Delsup{j-1}$
in \lex order, and indeed this coincides with computing the 
$-+$ rule on all of  $\xDel$ in \lex order because the other 
segments in  $\Delsup{k} $
only contribute blanks.
Using remark \ref{rem-cosoc} and the fact from the proof of lemma
\ref{lemma-ej} that 
 $\et{j} \Mg{ \Delsupj \cup \Delsup{j-1}}$ is the only composition
factor of 
 $\e{j} \Mg{ \Delsupj \cup \Delsup{j-1}}$
with $\epsj{j} = \epsj{j}(\Md) -1$, we get  
$$\et{j} \Md = \cosoc \Ind A \bx( \et{j} \Mg{ \Delsupj \cup \Delsup{j-1}})
\bx B = \cosoc \Ind (\E{j}\Del) = \Mg{ \E{j} \Del}.$$
   \end{proof}

   
The proof of theorem \ref{thm-Ehat},
that $\ethati \Md = \Mg{ \Ehati \Del}$,
is analogous to the proof of theorem \ref{thm-Ei},
and so will not be included.

As a result, we also get corollary \ref{cor-cycl}, which
uses rule \ref{defn-Ehat} to determine which $\Md$ are in
$\Rep \Hlamn$.
In the next section we will study the case $\lambda = \Lambda_i$.
\section{Multipartitions} 
\label{sec-proofpartition}

In this section, we will combine remark \ref{rem-cosoc} with
lemma \ref{lemma-triple} to build up $\Md$ from irreducible
$H_{n_i}^{\Lambda_i}$-modules.

We provide a dictionary between multisegments and certain colored
multipartitions, ``Kleshchev multipartitions''. 
The reverse map is obvious, but the forward one
is more subtle.
\omitt{because it is not surjective.}
  Only certain colored
multipartitions can arise.
In section \ref{sec-crystals} we explain those colored multipartitions
correspond to
nodes in a connected component of a tensor product of level 1 crystals.

We can use theorem \ref{thm-Ehat} to identify 
for which $\Del$
is $\Md \in \RHlam$.

\begin{thm}
\label{thm-hfin}
 The $\Del$ such that $\Md \in \RLami$ are 
of the form
$\Del = \{ \del i{b_1}, \del{i-1}{b_2}, \ldots,$  $ \del{i-k+1}{b_{k}} \}$,
where $b_1 > b_2 > \cdots > b_{k}$, or $\Del = \emptyset$.
\end{thm}
In other words, we can associate to $\Del$ the $i$-colored 
\omitt{(by $i$)}
 partition $\mu(\Del )= (b_1 -i+1, b_2 -i+2, \ldots, b_{k} -i +k)$,
pictured as 
$$ \tableau{i& 
{\scriptstyle i+1}
& \cdots&\mbox{}& \cdots& b_1 \\
	{\scriptstyle i-1}& i& \cdots& b_2\\
	\vdots& \mbox{}&\mbox{}\\
{\scriptscriptstyle i\!-\!k\!+\!1}& \cdots & b_{k}}.$$
In this case, $\Del(\mu, i) = \Del$.

\begin{proof}
To be in $\RLami$ means that $\epshati(\Md) \le 1$ and
$\epshat{j}(\Md)  = 0$ if $j \neq i$.
Rule \ref{defn-Ehat} for computing $\epshat{j}(\Md)$ shows
that if the segment $\del jb$ occurs in $\Del$ and $j \neq i$, then
it must be preceded in \std order by $\del{j+1}{b'}$, which forces
$b < b'$.
Further, if $\Del \neq \emptyset$ then some segment $\del i{b_1}$ must occur.
Only the $\Del$ listed above satisfy these requirements.
\end{proof}

It is no surprise that the multisegments above correspond to
partitions. We note that for generic $q$, $H_n^{\Lambda_0} \isom R S_n$,
the group algebra of the symmetric group.
The irreducible modules of $S_n$ are parameterized by partitions.
The irreducible module  parameterized by 
$\mu$ is called the Specht module $S^\mu$ and
it is isomorphic to $\Mg{\Del(\mu,0)}$.
(The modules for $H_n^{\Lambda_i}$ differ only in that $X_1$ acts by the
scalar $\qi$.)

Recall we had defined
$\Del(\mu, i) = \{ \del{i}{i+\mu_1-1}, \ldots, \del{i - k+1}{i-k +\mu_{k}} \}$,
where $k$ is the length of $\mu$.
Hence, 
$\mu\left( \Del(\mu,i) \right) = \mu$ undoes this operation.
Theorem \ref{thm-part1} explains how to associate
a multipartition to an arbitrary multisegment.

Also recall that $\Nmui = \Mg{\Del(\mu, i)}$,
with the convention that $\Nm{\0}{i} = 0$.


\begin{lemma}
\label{lemma-box}
\begin{enumerate}
\item For all $j$, $\epsj{j}(\Nmui) \le 1$. 
\item
$\et{j} \Nmui = \Nmb{\E{j}(\mu,\Lambda_i)}$,
where, as in Rule \ref{defn-p}, $\E{j}(\mu,\Lambda_i)$
denotes the removal of a removable $j$ box from the diagram
of $\mu$ when colored by $i$,
or denotes $\0$ if no such removable $j$ box exists.
\end{enumerate}
\end{lemma}
\begin{proof}
In the multisegment $\Del = \Del(\mu, i) =
\{ \del{i}{b_1}, \ldots, \del{i - k  +1}{b_k} \}$ with $b_1 > \cdots > b_k$,
observe  each $b_t$ is distinct.  Thus in our rule for computing
$\epsj{j}(\Md)$, there is at most one $-$ symbol.  This shows part 1
(using theorem \ref{thm-Ei}).

In fact, if $\epsj{j} ( \Nmui) = 1$, it means we must see  $-$ but
not  $-+$ in calculating as in rule \ref{defn-Ei}.
This corresponds to the $i$-colored diagram for $\mu$ having a removable
$j$-box.  The $-$ means some row of $\mu$ ends in a $j$-box.
But the absence of $+$ means the rows below it must not end in a $(j-1)$-box. 
That means the $j$-box has no box directly below it (else it would be filled
with $j-1$). Then it is clearly removable,
and $\Del(\mu \setminus \jboxj, i) = \E{j} \Del$.
The converse also holds, since if $\jboxj$ is removable,
the row beneath it cannot end in $j-1$ but nor can any of the
rows beneath it, since their fillings strictly decrease.
Thus $\et{j} \Nmui = \et{j} \Md = \Mg{\E{j} \Del}
=\Mg{\Del(\mu \setminus j,i)} = \Nmb{\E{j}(\mu,\Lambda_i)}$.


\end{proof}

\begin{thm}
\label{thm-two}
\begin{enumerate}
\item
Let $\Md$ be an $H_n^{\Lambda_i + \Lambda_\h}$-module.
Assume $i \ge \h$.
Then there exist partitions $\mu$ and $\nu$ such that
$$\Md = \cosoc \Ind \Nmui \bx \Nm{\nu}{\h}.$$
Further, 
 $\Del = \Del(\mu, i) \cup \Del(\nu, \h)$ and 
$b_{i-\h+x} \le c_x$, where 
$\Del(\mu, i) = \{ \del i{b_1}, \ldots, \del{i-m+1}{b_m} \}$ and
$\Del(\nu, \h) = \{ \del \h{c_1}, \ldots, \del{\h- k+1}{c_k} \}$.
\item
Given $\mu, \nu$ such that $\mu_{i-\h+x} \le \nu_{x}$ for $x \ge 1$,
$\cosoc \Ind \Nmui \bx \Nm{\nu}{\h}$ is irreducible and equal to
$\Md$ where $\Del = \Del(\mu, i) \cup \Del(\nu, \h)$.
\end{enumerate}
\end{thm}
\begin{proof}
Similar to lemma \ref{lemma-box}, $\Del$ must consist
of two consecutive (with respect to the first ``coordinate'' of a segment)
decreasing subsequences.  Applying the sliding lemma
will let us sort the segments comprising the subsequences  into two pieces corresponding to $\mu$ and $\nu$.

Write $\Del = \{ \del{a_1}{b_1}, \ldots, \del{a_t}{b_t} \}$ in \lex order.
We must have $\epshat{k} (\Md) = 0$ if $k \neq i, \h$, $\epshati(\Md) \le 1$
and $\epshat{\h}(\Md) \le 1$ (or $\epshati(\Md) \le 2$ if $i = \h$).
Theorem \ref{thm-Ehat} shows that $\{ a_1, \ldots, a_t \}$ is the union of
two intervals $\{ i, i-1, \ldots, i-s+1\}$ and $\{\h, \h-1, \ldots, \h -t+s+1\}$
where their corresponding $b$'s are also strictly decreasing.

If these intervals are disjoint, then the concatenation
of the intervals is in \lex order and it is clear
\begin{eqnarray*}
\Md &=& \cosoc \Ind \xDel \\
	&=& \cosoc \Ind \del i{b_1} \bx \del{i-1}{b_2} \bx \cdots \bx \del{i-s+1}{b_s} \bx
\del{\h}{b_{s+1}} \bx \cdots \bx \del{\h-t+s+1}{b_t} \\
	&=& \cosoc \Ind_{s,t}^n (\cosoc \Ind \del i{b_1} \bx \cdots \bx \del{i-s+1}{b_s}) \bx
(\cosoc \Ind \del{\h}{b_{s+1}} \bx \cdots \bx \del{\h-t+s+1}{b_t}) \\
	&=& \cosoc \Ind_{s,t}^n \Nmui \bx \Nm{\nu}{\h},
\end{eqnarray*}
where $\mu = (b_1 -i+1, \ldots, b_s -i+s)$ and 
$\nu = (b_{s+1} - \h +1, \ldots, b_t - \h +t -s)$.
By construction $\Del = \Del(\mu, i) \cup \Del(\nu, \h)$.

However, there may be overlap between the intervals.
In that case, we partition $\Del = \Del_1 \cup \Del_2$ as follows.
If $a_m \in \{ i, \ldots, i-s+1\} \setminus \{ \h, \ldots, \h-t+s+1\}$ then
put $\del{a_m}{b_m}$ in $\Del_1$.
If $a_m \in \{ \h, \ldots, \h-t+s+1\} \setminus \{ i, \ldots, i-s+1\} $ then
put $\del{a_m}{b_m}$ in $\Del_2$.
For $a \in \{ i, \ldots, i-s+1\} \cap \{ \h, \ldots, \h-t+s+1\}$, 
there exist $c \ge b$ such  that both
of $\del ab$ and $\del ac$ are in the multisegment $\Del$.
Put the shorter $\del ab$ in $\Del_1$ and the longer $\del ac$ in $\Del_2$.
Then $\mu = \mu(\Del_1)$ and $\nu = \mu(\Del_2)$ will be partitions and
$\Del_1 = \Del(\mu, i),$ $\Del_2 = \Del(\nu, \h)$.

Define $\Gamma_1^a = \{ \del \alpha \beta \in \Del_1 \mid \alpha > a\}$
and $\Gamma_2^a = \{ \del \alpha \beta \in \Del_2 \mid \alpha > a\}$.
Suppose $b \le c$.
Then it follows from part 2 of the sliding lemma \ref{lemma-triple} that
\begin{align*}
\cosoc &\Ind \Mg{\Gamma_1^a} \bx \Mg{\Gamma_2^a} \bx \del ac \bx \del ab \\
&= \cosoc \Ind \Mg{\Gamma_1^a} \bx \Mg{\Gamma_2^a \cup \{\del ac\}} 
	\bx \del ab \\
&=  \cosoc \Ind \Mg{\Gamma_1^a} \bx  \del ab \bx \Mg{\Gamma_2^a \cup 
	\{ \del ac\} }  \\
&=  \cosoc \Ind \Mg{\Gamma_1^a \cup \{ \del ab\}} \bx  \Mg{\Gamma_2^a \cup 
	\{\del ac \} }
&=  \cosoc \Ind \Mg{\Gamma_1^{a-1} } \bx  \Mg{\Gamma_2^{a-1}}.
\end{align*}
If $a$ is such that both $\del ab, \del ac \in \Del$, then $\Gamma_1^{a-1} = 
\Gamma_1^a \cup \{ \del ab \}$ and 
$\Gamma_2^{a-1} =  \Gamma_2^a \cup \{ \del ac \}$. 
Continuing as above, starting with $\xDel$ in \lex order, yields
$\Md = \cosoc \Ind \Mg{\Del_1} \bx \Mg{\Del_2} = \cosoc \Ind \Nmui \bx \Nm{\nu}{\h}$.

Conversely, given any $\mu, \nu$ with the property that 
$\mu_{i-\h+x} \le \nu_x$ for all $x \ge 1$ (with the convention
$\mu_x = 0 $ if $x > \ell(\mu)$),
the above argument shows $\cosoc \Ind \Nmui \bx \Nm{\nu}{\h}$ is irreducible
and equal to $\Md$ for $\Del = \Del(\mu, i) \cup \Del(\nu, \h)$.

Observe that the condition that $\mu$ and $\nu$ satisfy are equivalent
to the multipartiton  $(\mu, \nu)$
colored by $ \Lambda_i + \Lambda_h$ being a Kleshchev multipartition.
\end{proof}

Observe that $n(\Del) = |\mu| + |\nu|$.  Also it is possible that $\mu$ or $\nu$ 
be the empty partition.

Theorem \ref{thm-part1} extends the same argument to $\lambda$ of any level.
We repeat its statement here for convenience.
\omitt{repeat?}
\begin{thmA}[theorem \ref{thm-part1}]
\label{thm-part12}
Let $\Md$ be an $\Hlamn$-module where $\lambda = \Lambda_{i_1} +
\Lambda_{i_2} + \cdots + \Lambda_{i_r}$, $i_1 \ge \cdots \ge i_r$.
Then there exists
a Kleshchev multipartition $\mubarik$
\omitt{
$\underline{\mu} = (\mu^{(1)}, \ldots, \mu^{(r)})$ satisfying 
\begin{gather}
\label{eq-connected}
\mu_{i_t - i_{t+1} + x}^{(t)} \le \mu_x^{(t+1)} \quad
\text{ for all $x \ge 1$, $1 \le t \le r-1$}
\end{gather}
such that $\Md = \cosoc \Ind \Nm{\mu^{(1)}}{i_1} \bx \cdots \bx \Nm{\mu^{(r)}}{i_r}$
and $\Del = \Del(\mu^{(1)}, i_1) \cup \cdots \cup \Del(\mu^{(r)}, i_r)$.
}
such that $\Md = \Nmulam$ and $\Del = \Del(\mubarlam)$.

\omitt{
Conversely, if $\mubari$ satisfies \eqref{eq-connected}, then
$\cosoc \Ind \Nm{\mu^{(1)}}{i_1} \bx \cdots \bx \Nm{\mu^{(r)}}{i_r}$
is irreducible and isomorphic to $\Md$ for $\Del= \Del(\mu^{(1)}, i_1)
\cup \cdots \cup \Del(\mu^{(r)}, i_r)$.
}
Conversely, if $\mubarik$ 
is a Kleshchev multipartition, then 
$\Nmulam = \cosoc \Ind \Nm{\mu^{(1)}}{i_1} \bx \cdots \bx \Nm{\mu^{(r)}}{i_r}$
is irreducible and isomorphic to $\Md$ for $\Del= \Del(\mubarlam)$.

 \end{thmA}

If $\mubari$  
is not a Kleshchev multipartition, then
$\cosoc \Ind \Nm{\mu^{(1)}}{i_1} \bx \cdots \bx \Nm{\mu^{(r)}}{i_r}$
need not be  irreducible.
Furthermore, Theorem \ref{thm-part2} (repeated below) need not
hold for arbitrary colored multipartitons.

The reader may recognize the condition
\begin{gather*}
\mu_{i_t - i_{t+1} + x}^{(t)} \le \mu_x^{(t+1)}
\end{gather*}
as the condition that the nodes in the connected component
of highest weight $\lambda$ satisfy in the tensor of crystal graphs
$\B(\Lambda_{i_r})\otimes \cdots \B(\Lambda_{i_1})$.
See section \ref{sec-crystals} and corollary \ref{cor-part4tensor}
for a description of these theorems
in the language of crystal graphs.

We will now justify that the $\pm$ rule \ref{defn-Ei} given
for $\Ei \Del$  is compatible with the one for multipartitions.

\begin{thmA}[theorem \ref{thm-part2}]
\label{thm-part22}
Given a multisegment $\Del$ and a Kleshchev multipartition
$\mubari$ such that $\Del = \Del(\mubari)$, then
$\E{j} \Del = \Del( \E{j} \mubari)$.
In other words, $\et{j} \Nmb{\mubari} = \Nmb{ \E{j} \mubari}$.
\end{thmA}

\begin{proof}
When we compute $\E{j} \Del$, we require that $\xDel$ be
in \lex order. However, in  constructing  $\mubari$ in theorem
\ref{thm-part1}, 
that order is changed by repeated application
of the sliding lemma. 

We'll consider all possible re-orderings that occur and show
they do not change the uncanceled $+$ and $-$ symbols.
There are 3 cases.

Case 1. In $\xDel$ we see  $\del a{j-1}$ followed by $\del aj$.
Then these two segments remain in that relative order.
(Recall in partitioning $\Del$ into two pieces that
our rule would put $\del a{j-1}$ in $\Del_1$ and
$\del aj$ in $\Del_2$.)

Case 2.
In $\xDel$ we see  $\del b{j-1}$ followed by $\del aj$ and $b > a$.
Again, these two segments remain in that relative order.

Case 3.
In $\xDel$ we see  $\del b{j}$ followed by $\del a{j-1}$ and $b > a$.
If these two segments do switch position, it would be by applying
the sliding lemma---that $\Ind \Mg{\Gamma} \bx \del a{j-1} \isom
\Ind \del a{j-1} \bx \Mg{\Gamma}$ where  $\Gamma \ni \del bj$ and the
last segment in $\Gamma$ is $\del ac$ for some $j > c \ge j-1$.
Thus, $c = j-1$.  Before sliding the pattern of symbols is $-++$,
whereas afterwards it is $+-+$.  However, all adjacent
$-+$ get cancelled in either rule for $\E{j}$ and so we have not
changed the uncanceled symbols (nor the segments they are attached to)
at all; both configurations
contribute just a $+$ attached to the $\del a{j-1}$.

Applying the $\pm$ rule to compute $\et{j} \Md = \Mg{ \E{j} \Del}$
is  thus equivalent to that of computing $\et{j} \Nmb{\mubari}
= \Nmb{\E{j} \mubari}$.
\end{proof}

Above we have built up the partitions ``row by row'', with
each row corresponding to a segment.  We could have chosen
to build them ``column by column'' which would have meant using
Steinberg modules in place of trivial modules.
These are also one dimensional, but have character
$(q^j q^{j-1} \cdots q^i)$,
and all $T_k +1$ vanish on them.
In that case we would have a similar $+-$ rule that would correspond to looking
for addable/removable boxes at the ends of columns (not rows).

\section{Discussion: crystals}
\label{sec-crystals}
Now we can put the theorems of section \ref{sec-main}
into the language of Kashiwara's  crystal graphs. 
We refer the reader to \cite{K} for the definitions and properties
of crystal graphs.
A good reference in this context is  \cite{G},
where he also includes when $q$ is a root of unity, in which
case one must be more careful.

Fix $\lambda = \Lambda_{i_1} + \cdots + \Lambda_{i_r}$
with $i_1 \le i_2 \le \cdots \le i_r$, a weight of the
Lie algebra $\gl$.
Let $L(\lambda)$ denote the irreducible integrable  representation
of $\gl$
with highest weight $\lambda$ and let $\B(\lambda)$ denote
its crystal graph.

Let $\B(\infty)$ denote the crystal graph associated to 
$U(\eta^-)$. 
 We can think of each node of  $\B(\infty)$  as labeled
by a multisegment $\Del$.
(Observe another interpretation of a multisegment is
as a sum of positive roots of $\gl$.
A segment $\del ij$ 
corresponds to the positive root
$\alpha_{ij} = \alpha_i + \cdots + \alpha_j$,
and the multisegment $\Del$  corresponds to their sum.)

The following two theorems are essentially Grojnowski's Theorem 14.3.
\begin{thm}\cite{G}
\label{thm-graph}
Let 
$\B_\lambda$
be the graph
 whose nodes are the irreducible 
$\Hlamn$-modules for $n \ge 0$ with edges given
by $\eti M \xrightarrow{i} M$ if $\eti M \neq 0$.
Then $\B_\lambda = \B(\lambda)$.
\end{thm}
\begin{thm} 
\label{thm-graphinfty}
Let 
$\B_\aff$
be the graph
whose nodes are the irreducible 
modules in $\Rq$ with edges given
by $\eti M \xrightarrow{i} M$ if $\eti M \neq 0$.
Then $\B_\aff = \B(\infty)$.
\end{thm}



In this language, 
Theorems \ref{thm-Ei}, \ref{thm-part1}, and \ref{thm-part2} become
\begin{cor}
$\B_{BZ} =\B_\aff$ 
with crystal isomorphism $\Del \mapsto \Md$.

$\B'_\lambda = \B_\lambda$ 
with crystal isomorphism 
$\mubarlam \mapsto  \Nmulam$.

Furthermore,
$\mubarlam \mapsto \Del(\mubarlam)$
coincides with the inclusion of $\B_\lambda$ into $\B_\aff$.
\end{cor}

We can view $\B_\lambda$ as a subgraph of $\B_\aff$
(or $\B(\lambda) \subseteq \B(\infty)$).
Indeed the irreducible $\Hlamn$-modules are just a subset
of the irreducibles of $\RqHaffn$.
Rule \ref{defn-Ehat} for $\epshati (\Md)$ tells us exactly  which subgraph
 $\B_\lambda$ is, i.e.~ which multisegments comprise it,
by requiring $\epshati(\Md)$ be less than or equal to the multiplicity
of $\Lambda_i$ in $\lambda$.

On the other hand, $\B(\lambda) \inj \B(\Lambda_{i_r}) \otimes
\cdots \otimes \B(\Lambda_{i_1})$.
The nodes in the graph  $\B(\Lambda_i)$ correspond to
{\it all\/} partitions $\mu$, colored by $i$.
(In the case $i=0$, $L(\Lambda_0)$ is the basic representation
and the crystal graph
$\B(\Lambda_0)$ can be identified with Young's lattice
of partitions.  For other $i \in \Z$ simply shift the edge
labels by $+i$.)
Thus, in the tensor product $ \B(\Lambda_{i_r}) \otimes
\cdots \otimes \B(\Lambda_{i_1})$, it is natural to label nodes
by colored multipartitions $\mubari$.
The nodes run over all $\lambda$-colored multipartitions, but
this graph is not connected.

The connected component in $ \B(\Lambda_{i_r}) \otimes
\cdots \otimes \B(\Lambda_{i_1})$ of $\emptyset \otimes \emptyset
\otimes \cdots \otimes \emptyset$,
the unique node with
weight $\lambda$, defines $\B(\lambda)$ as a subgraph.
These realizations of $\B(\lambda)$  as a subgraph of two different
crystals
gives us a map  from multisegments (those with $\Ehati$ bounded
by $\lambda$)  to colored multipartitions
\begin{gather*}
\Del \mapsto \mu^*(\Del)
\end{gather*}
which respects the action of $\eti$.

As $r$ goes to $\infty$ the domain ranges over all multisegments.
However we never get all colored multipartitions in the image. 
Theorem \ref{thm-part1} explains which ones we do get, and
corollary \ref{cor-part4tensor} 
tells us how to interpret those
in terms of a tensor product of crystal graphs.

There is a small technical point here which we must address.
The map above  is the ``reverse'' of the map described in 
theorem \ref{thm-part1}.
  One can simply think of
$(\mubari)^*$ as $\mubari$ read from right to left.
%
The conventional definition for tensoring crystals and the
conventional definition for partitions requires that we introduce
a  reversal in order to be consistent.
We choose to modify the convention for tensoring crystals
instead of the convention of reading a partition from left to
right.
(One way we could have changed the definitions here
to modify  instead the convention for partitions would be
if we had taken segments
in increasing order 
and  instead of taking $\cosoc \Ind$ everywhere we had taken
$\soc \Ind$ (or $\cosoc \Indhat$), by proposition \ref{prop-soccosoc}.
See remark \ref{rem-hat} for further comments.)

This reversal motivates the following definition, which reverses
the order of tensoring crystal graphs.

 Let $\B_1 \ot \B_2 =\B_2 \otimes \B_1$, and so $\eti$ acts via
\begin{gather}
\label{eqn-crystal2}
\eti (b_1 \ot  b_2) = \begin{cases}
	 \eti b_1 \ot  b_2 
	  &\text{ if $\epsi(b_1) > \varphi_i(b_2)  $}\\ 
	b_1 \ot  \eti b_2 
	  &\text{ if $\epsi(b_1) \le \varphi_i(b_2).$}
\end{cases}
\end{gather}

%
The usual convention for the action of $\eti$ on 
 $\B_1 \otimes \B_2$ is:
\begin{gather}
\label{eqn-crystal1hat}
\eti (b_1 \otimes b_2) = \begin{cases}
	 \eti b_1 \otimes b_2 & \text{ if $\varphi_i(b_1) \ge \epsi(b_2)$}\\
	b_1 \otimes \eti b_2 & \text{ if $\varphi_i(b_1) < \epsi(b_2)$.}
\end{cases}
\end{gather}
Either way is compatible with the following picture
(in which all edges $\to$
are assumed labelled by $i$).

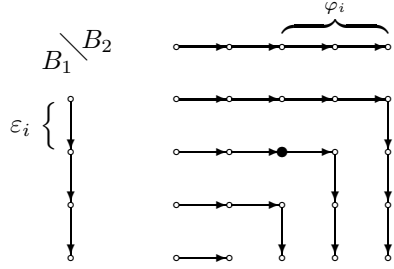
\begin{center}
\begin{picture}(100,110)(0,-10)
\multiput(40,0)(20,0){5}{\circle{2}}
\multiput(40,20)(20,0){5}{\circle{2}}
\multiput(40,40)(20,0){5}{\circle{2}}
\multiput(40,60)(20,0){5}{\circle{2}}
\multiput(40,80)(20,0){5}{\circle{2}}
\multiput(41,80)(20,0){4}{\vector(1,0){18}}
\multiput(41,60)(20,0){4}{\vector(1,0){18}}
\multiput(41,40)(20,0){3}{\vector(1,0){18}}
\multiput(41,20)(20,0){2}{\vector(1,0){18}}
\put(41,0){\vector(1,0){18}}
\multiput(0,0)(0,20){4}{\circle{2}}
\multiput(0,59)(0,-20){3}{\vector(0,-1){18}}
\multiput(120,59)(0,-20){3}{\vector(0,-1){18}}
\multiput(100,39)(0,-20){2}{\vector(0,-1){18}}
\put(80,19){\vector(0,-1){18}}
\put(80,40){\circle*{4}}
\put(80,81){\makebox(40, 10)[s]{$\overbrace{\makebox(40,10){}}^{\varphi_i}$}}
\put(-10,40){\makebox(5, 20)[r]{$\varepsilon_i $ $\Bigl\{ $}}
\put(-4,85){\line(1,-1){10}}
\put(8,78){\makebox(8, 8)[r]{$B_2$} }
\put(-7,70){\makebox(8, 8)[r]{$B_1$} }

\end{picture}
\end{center}

This is usually a picture
\omitt{of $B \otimes A$, whereas here we think of it as $A \ot B$.
And so the solid dot is $a \ot b \in A \ot B$,
or $b \otimes  a \in B \otimes  A$,}
of $B_2 \otimes B_1$, whereas here we think of it as $B_1 \ot B_2$.
And so the solid dot is $b_1 \ot b_2 \in B_1 \ot B_2$,
or $b_2 \otimes  b_1 \in B_2 \otimes  B_1$,
and the rule above tells you if an arrow leading to it
approaches from the top or the left.

In general, the crystals $B_1 \ot B_2$ and $B_1 \otimes B_2$ are
not isomorphic; but in this setting, since the crystal graphs
correspond to integrable highest weight representations of $\gl$,
they are isomorphic. Reversing the order here is preferable because
then the node labelled by $\Nmulam$ in $\B_\lambda$ will correspond
to the node 
$\mu^{(1)} \ot  \mu^{(2)} \ot  \cdots \ot  \mu^{(r)}$
in $B(\Lambda_{i_1}) \ot  B(\Lambda_{i_2}) \ot 
\cdots  \ot  B(\Lambda_{i_r})$ (read now from left to right).
Those readers more familiar with partitions may prefer this order,
while those more familiar with crystal graphs may find it annoying.

\omitt{This is compatible with \std order.}%
Compare the definition for $\B_1 \ot \B_2$
 to Theorems \ref{thm-Ei} and \ref{thm-part2},
which imply that  for very special $A$ and $B$, 
\begin{gather}
\label{eqn-crystal3}
\eti (\cosoc \Ind A \bx B) = \begin{cases}
	\cosoc (\Ind \eti A \bx B) 
 &\text{ if $\epsi(A) > \varphi_i(B)$ }\\ 
	\cosoc (\Ind A \bx \eti B) 
 &\text{ if $\epsi(A) \le \varphi_i(B)$,  }
\end{cases}
\end{gather}
and also that the cosocles above are irreducible.

\medskip

We will now explain \eqref{eqn-crystal3} in more detail and in the process
prove corollary \ref{cor-part4tensor}.

%

First we will define $\vpi$.
We refer the reader to \cite{G} to learn more about 
$\varphi_i$, which is more subtle than $\epsi$.
It depends on 
the crystal operator $\fti$, mentioned only briefly here
and used to simplify the proof of corollary \ref{cor-alpha} below.

The crystal operator $\fti$ satisfies 
$$\eti M = N \iff M = \fti N.$$
Module-theoretically, $\fti N  = \cosoc \Ind_{n,1}^{n+1} N \bx \qi
 = \cosoc \Ind N \bx \del ii$.
This gives a recipe of how to construct an irreducible
$\Haffn$-module $M$ using the crystal graph $\B_\aff$,
 by following any path on the crystal from
the node corresponding to $M$ back to the root,
one step at a time.
In other words, given
$\1 = \et{j_n}  \cdots \et{j_1} M$, we have
$M = \ft{j_1}  \cdots \ft{j_n} \1$.
As in remark \ref{rem-crystal},
\cite{BZ,Z} already gave a distinguished such path.
For more on the properties of  $\fti$ see \cite{G}.

For $\B_\lambda,$
 $\varphi_i(M) = \max\{m \ge 0 \mid \fti^m M \in \Rep \Hlamn\}$.
For $\B_\aff$ it is slightly different.
However, in this context, $\varphi_i$ is computed 
by counting uncanceled $+$ symbols as in Rule \ref{defn-Ei}
or \ref{defn-p}, just as $\epsi$  
counts uncanceled $-$ symbols.

If we decree that  $\varphi_i$ counts uncancelled $+$ symbols,
then we can see
when Theorems \ref{thm-Ei} and \ref{thm-part2} imply
 \eqref{eqn-crystal3}, as follows.
 
First, the $A$ and $B$ we consider are of the form
$A= \Mg{\Gamma}$ and $B= \Mg{\Gamma'}$, where now
$\Gamma$ is an initial segment of $\Del$ and 
$\Gamma'$ is a final segment, with respect to 
\lex order, and $\Gamma \cup \Gamma' = \Del$.
(Note the different use of the words initial  ``segment'' here is
as initial subword  with respect to \lex order.)
Then $\Ind \Del = \Ind \Gamma \bx \Gamma'$.
By transitivity of induction and remark \ref{rem-cosoc},
$\Md = \cosoc \Ind \Mg{\Gamma} \bx \Mg{\Gamma'}$.
Theorem \ref{thm-Ei} says that the $\pm$ rule works
for all three of $\Md, \Mg{\Gamma}$, and $\Mg{\Gamma'}$.
The $\pm$ word for $\Del$ is the concatenation of
the $\pm$ words for $\Gamma$ and $\Gamma'$.
$$++ \cdots + \underbrace{- - \cdots -}_{\epsi} \quad 
\underbrace{++ \cdots +}_{\vpi} - - \cdots -$$
If one then performs  the requisite cancelling of $-+$ pairs,
the leftmost uncanceled $-$
for $\Gamma$ (which has a total of
$\epsi(\Mg{\Gamma})$ as yet  uncanceled $-$)
 will remain such in $\Del$ only if it is not followed
by more than $\epsi(\Mg{\Gamma})$ many as yet uncancelled  $+$, i.e.~
if   $\epsi(\Mg{\Gamma}) > \varphi_i(\Mg{\Gamma'})$.
But in the case $\epsi(\Mg{\Gamma}) \le \varphi_i(\Mg{\Gamma'})$,
all $-$ from $\Gamma$ will get cancelled by subsequent $+$,
and it is now the leftmost uncanceled $-$ from $\Mg{\Gamma'}$ that
becomes the leftmost uncanceled $-$ for ${\Del}$.
On the one hand, we are describing the rule for $\eti \Md$, and
on the other hand, we are describing \eqref{eqn-crystal3}.
Notice at each stage, that all cosocles we take are irreducible
(since the conditions of \eqref{eqn-crystal3} ensure the concatenation
$\Ei(\Gamma) \cup \Gamma'$ or $\Gamma \cup \Ei(\Gamma')$ is in \lex order).

Next, the argument for $A$ and $B$ of the form $\Nmulam$ is similar;
but again, the multipartitions must be  initial and final
segments of a given multipartition,
and all three must be Kleshchev multipartitions, i.e.~ colored
multipartitions
 satisfying condition \eqref{eq-connected1}.

Consequently, we also get the following corollary,
which describes how to build irreducible $\Hlamn$-modules
out of $H_a^\alpha$ and $H_b^\beta$-modules when $\alpha + \beta = \lambda$.
\begin{cor}
\label{cor-alpha}
Suppose $\alpha$ and $\beta$ are weights with $\alpha + \beta = \lambda$
and $\alpha = \sum_{k=1}^a \Lambda_{i_k}$,
$\beta = \sum_{k=a+1}^r \Lambda_{i_k}$, where $i_1 \le i_2 \le \cdots \le i_r$.
Suppose
$\mubarik \ot \nubarik \in \B(\alpha) \ot \B(\beta)$
is in the connected component of the unique node with  weight $\lambda$.
Write
$\concat{\mubarik}{\nubarik}
= (\mu^{(1)}, \ldots, \mu^{(a)}, \nu^{(1)}, \ldots, \nu^{(r-a)})$.
Then the module
$$ \Nmb{ \concat{\mubarik}{\nubarik} , \lambda} :=
\cosoc \Ind \Nmb{\mubaralpha} \bx \Nmb{\nubarbeta}$$
is irreducible
and
$$\eti (\Nmb{ \concat{\mubarik}{\nubarik} , \lambda}) =
\begin{cases}
	 \Nmb{\concat{\Ei(\mubarik)}{\nubarik}, \lambda}
&\text{ if $\epsi(\Nmb{\mubaralpha}) > \varphi_i(\Nmb{\nubarbeta}) $}\\ 
	\Nmb{\concat{\mubarik }{ \Ei(\nubarik)}, \lambda}	
&\text{ if $\epsi(\Nmb{\mubaralpha}) \le \varphi_i(\Nmb{\nubarbeta}) $}\\ 
\end{cases}$$
\end{cor}
\begin{proof}
The proof follows from 
the discussion after equation \eqref{eqn-crystal3}.
If we take the $\pm$ word corresponding to $\concat{\mubarik}{\nubarik},\lambda$
but then divide that word according to how $\alpha + \beta = \lambda$,
the rule for computing $\Ei$ coincides with the rule for
computing $\Ei( \mubarik \ot \nubarik)$ in $\B(\alpha) \ot \B(\beta)$.

We know by Theorem \ref{thm-part2} that some sequence
of $\eti$ will bring  
$\Nmb{ \concat{\mubarik}{\nubarik} , \lambda} $ 
to $\1 = \Nm{(\emptyset,\ldots,\emptyset)}{\lambda}$
if $\concat{\mubarik}{\nubarik} $ satisfies condition \eqref{eq-connected1}.
The argument above shows the same sequence of $\eti$ will
bring $\mubarik \ot \nubarik$ to $\emptyset \ot \emptyset$.
Thus Kleshchev multipartitions correspond to nodes in the connected
component of $\emptyset \ot  \emptyset$.

To see that every node in that connected component corresponds to
some Kleshchev multipartition, we need 
to verify that 
if $\eti (\mubarik'  \ot  \nubarik') = \mubarik  \ot  \nubarik$ 
and  $\concat{\mubarik}{\nubarik},\lambda $ is Kleshchev,
then so is
$\concat{\mubarik'}{\nubarik'},\lambda $.

This is less cumbersome to state if we make use of the
crystal operator $\fti$.
 We define $\Fi (\mubarlam)$ to
add an $i$-box to the partition corresponding to the
rightmost uncanceled $+$, if it exists, and otherwise it is $\0$.
Thus
$\Ei(\mubar{\mu'}{\lambda})  = \mubarlam \iff \mubar{\mu'}{\lambda}=
\Fi (\mubarlam).$
The argument about $\pm$ words above already shows that
$\Fi (\concat{\mubarik}{\nubarik},\lambda) =
\concat{\mubarik'}{\nubarik'},\lambda $ if 
$\fti(\mubarik \ot \nubarik) = \mubarik' \ot \nubarik'$,
where
\begin{gather}
\label{eqn-crystalf}
\fti (b_1 \ot  b_2) = \begin{cases}
	 \fti b_1 \ot  b_2 
	  &\text{ if $\epsi(b_1) \ge \varphi_i(b_2)  $}\\ 
	b_1 \ot  \fti b_2 
	  &\text{ if $\epsi(b_1) < \varphi_i(b_2).$}
\end{cases}
\end{gather}
This definition is also compatible with the picture 
below \eqref{eqn-crystal1hat}, and ``undoes'' $\eti$.

What we need to show is that if 
$\concat{\mubarik}{\nubarik},\lambda$ is Kleshchev, then so is
$\Fi (\concat{\mubarik}{\nubarik},\lambda)$, if it is nonzero.
In fact, given the comment above,
all we need show
is that $\mubarlam$ is
Kleshchev $\iff$ $\Fi(\mubarlam)$ is.

It is quite tedious to check this directly.
 Theorem \ref{thm-part1} let's us show this indirectly.
The theorem maps a multisegment to a {\it Kleshchev\/} multipartition,
and Theorem \ref{thm-part2} says this map intertwines 
the operator $\Ei$. Hence it also intertwines $\Fi$ (except
possibly when $\Fi(\mubarlam) = \0$, but in this case there is nothing
to check).
Hence $\Fi(\mubarlam)$ must be Kleshchev since it comes from some
multisegment.
\end{proof}

We point out that when $A$ and $B$ are  arbitrary modules,
we needn't have $\cosoc \Ind A \bx B$ be irreducible, nor
have  \eqref{eqn-crystal3} hold.
Furthermore,
although $\B(\lambda)$ is a subgraph of
$\B(\alpha) \ot \B(\beta)$
if $\alpha + \beta = \lambda$, there is no reason to expect
\eqref{eqn-crystal3} should hold for arbitrary $A \in \Rep H_a^\alpha$,
$B \in \Rep H_b^\beta$.  We are only certain it holds
when $\alpha$ is an initial segment of $\lambda$ and 
$\beta$ a final segment  (in the sense that
$\lambda$ corresponds to $i_1 i_2 \cdots i_r$, with $i_1 \le \cdots \le i_r$)
and when $A$,$B$ corresponds to a node in the connected component
with highest weight $\lambda$.
We also point out that the ability to order segments or partitions, as
we do, relies on $q$ being generic.

The argument above for $\Del = \Gamma \cup \Gamma'$
(or $\lambda = \alpha + \beta$)
 also works 
if we chop a $\pm$ word for a Kleshchev 
multipartition 
into $r$ pieces.
(The case $r=2$ is above in corollary \ref{cor-alpha}.)
Then rule \ref{defn-p} for computing $\Ei(\mubarlam)$ is exactly the rule for 
computing  the crystal operator
$\et{i}$ on the element 
$\mu^{(1)} \ot  \mu^{(2)} \ot  \cdots \ot  \mu^{(r)}$
$\in B(\Lambda_{i_1}) \ot  B(\Lambda_{i_2}) \ot 
\cdots  \ot  B(\Lambda_{i_r})$.


The reason  we are in the connected component of $\emptyset \ot 
\cdots \ot \emptyset$
is the same as above.
This then sketches the proof of corollary \ref{cor-part4tensor},
that the inclusion 
$\mubarlam \mapsto \mubarlam$
of Kleshchev multipartitions 
into {\it all\/}  $\lambda$-colored multipartitions
coincides with the inclusion of 
$\B(\lambda)$ in
$ B(\Lambda_{i_1}) \ot  B(\Lambda_{i_2}) \ot \cdots  \ot  B(\Lambda_{i_r})$.


In fact, the ``true'' definition of Kleshchev multipartitions,
and one that extends to $q$ a root of unity,
is that they are colored multipartitions corresponding
 to nodes in that connected component.
For a condition parallel to \eqref{eq-connected1} see \cite{JMMO}.
\begin{rem}
\label{rem-hat}
We make a final comment about the right/left disparity encountered.
The usual tensoring of crystals {\it is\/} compatible with
computing $\Ehati$ and \std order.
Compare equation \eqref{eqn-crystal1hat} to
\begin{gather*}
\label{eqn-crystal2hat}
\ethati (\cosoc \Ind \Mg{\Gamma} \bx \Mg{\Gamma'}) = \begin{cases}
	\cosoc (\Ind \ethati \Mg{\Gamma} \bx \Mg{\Gamma'}) 
 &\text{ if $ \phati(\Mg{\Gamma}) \ge \epshati(\Mg{\Gamma'}) $ }\\ 
	\cosoc (\Ind \Mg{\Gamma} \bx \ethati \Mg{\Gamma'}) 
 &\text{ if $ \phati(\Mg{\Gamma}) < \epshati(\Mg{\Gamma'})$. }\\ 
\end{cases}
\end{gather*}
which happens in the case
$\Gamma$ is an initial segment of $\sDel$ and 
$\Gamma'$ is a final segment, with respect to 
\std  order, and $\Gamma \cup \Gamma' = \sDel$.

Just as $\epshati$ counted the number of uncanceled $-$ signs
when computing $\Ehati$ in Rule \ref{defn-Ehat},
$\phati$ counts the number of uncanceled $+$ signs.
\end{rem}



  \end{document}